\documentclass[12pt,british]{article}

\usepackage{babel}
\usepackage{amsmath,amsfonts,amssymb,amsthm}
\usepackage[latin1]{inputenc}
\usepackage[T1]{fontenc}
\usepackage{enumerate}
\usepackage{textcomp,url}
\usepackage{eucal}

\usepackage[noBBpl]{mathpazo}
\usepackage{mathabx}

\usepackage[matrix,arrow]{xy}

\makeatletter

\renewcommand{\p@enumii}{}

\def\@enum@{\list{\csname label\@enumctr\endcsname}%
           {\usecounter{\@enumctr}\def\makelabel##1{
\normalfont\ignorespaces\emph{{##1}~}}
\setlength{\labelsep}{3pt}
\setlength{\parsep}{0pt}
\setlength{\itemsep}{0pt}
\setlength{\leftmargin}{0pt}
\setlength{\labelwidth}{0pt}
\setlength{\listparindent}{\parindent}
\setlength{\itemsep}{0pt}
\setlength{\itemindent}{0pt}
\topsep=3pt plus 1pt minus 1 pt}}

\makeatother

\renewcommand{\epsilon}{\ensuremath{\varepsilon}}
\renewcommand{\phi}{\ensuremath{\varphi}}
\newcommand{\vide}{\ensuremath{\varnothing}}
\renewcommand{\to}{\ensuremath{\longrightarrow}}

\newcommand{\N}{\ensuremath{\mathbb N}}
\newcommand{\Z}{\ensuremath{\mathbb Z}}
\newcommand{\dt}{\ensuremath{\mathbb D}^{2}}
\newcommand{\St}[1][2]{\ensuremath{\mathbb S}^{#1}}
\newcommand{\FF}{\ensuremath{\mathbb F}}
\newcommand{\F}[1][n]{\ensuremath{\FF_{{#1}}}}
\newcommand{\rp}{\ensuremath{\mathbb{R}P^2}}

\newcommand{\sn}[1][n]{\ensuremath{S_{{#1}}}}
\newcommand{\an}[1][n]{\ensuremath{A_{{#1}}}}

\DeclareRobustCommand*{\up}[1]{\textsuperscript{#1}}
\renewcommand{\th}{\ensuremath{\up{th}}}
\newcommand{\ft}[1][n]{\ensuremath{\Delta_{#1}}}
\newcommand{\garside}[1][n]{\ensuremath{T_{#1}}}

\renewcommand{\ker}[1]{\ensuremath{\operatorname{\text{Ker}}\left({#1}\right)}}
\newcommand{\im}[1]{\ensuremath{\operatorname{\text{Im}}\left({#1}\right)}}
\newcommand{\aut}[1]{\ensuremath{\operatorname{\text{Aut}}\left({#1}\right)}}
\newcommand{\diff}[1]{\ensuremath{\operatorname{\text{Diff}^+}\!\!\left({#1}\right)}}
\newcommand{\id}{\ensuremath{\operatorname{\text{Id}}}}

\newcommand{\dih}[1]{\ensuremath{\operatorname{\text{D}}_{#1}}}
\newcommand{\dic}[1]{\ensuremath{\operatorname{\text{Dic}}_{#1}}}
\newcommand{\quat}[1][8]{\ensuremath{\mathcal{Q}_{#1}}}

\makeatletter
\def\@map#1#2[#3]{\mbox{$#1 \colon\thinspace #2 \to #3$}}
\def\map#1#2{\@ifnextchar [{\@map{#1}{#2}}{\@map{#1}{#2}[#2]}}
\makeatother

\newcommand{\brak}[1]{\ensuremath{\left\{ #1 \right\}}}
\newcommand{\ang}[1]{\ensuremath{\left\langle #1\right\rangle}}
\newcommand{\set}[2]{\ensuremath{\left\{#1 \,\mid\, #2\right\}}}
\newcommand{\setang}[2]{\ensuremath{\ang{#1 \,\mid\, #2}}}
\newcommand{\setangr}[2]{\ensuremath{\ang{#1 \,\left\lvert \, #2 \right.}}}

\newcommand{\setr}[2]{\ensuremath{\brak{#1 \,\left\lvert \, #2 \right.}}}
\newcommand{\setl}[2]{\ensuremath{\brak{\left. #1 \,\right\rvert \, #2}}}

\newtheoremstyle{theoremm}{}{}{\itshape}{}{\scshape}{.}{ }{}

\theoremstyle{theoremm}
\newtheorem{thm}{Theorem}

\newtheorem{prop}[thm]{Proposition}

\newtheoremstyle{remarkk}{}{}{}{}{\scshape}{.}{ }{}

\theoremstyle{remarkk}

\newtheorem{rem}[thm]{Remark}
\newtheorem{rems}[thm]{Remarks}

\newtheoremstyle{comment}{}{}{}{}{\bfseries}{:}{ }{}
\theoremstyle{comment}

\newcommand{\reth}[1]{Theorem~\protect\ref{th:#1}}

\newcommand{\repr}[1]{Proposition~\protect\ref{prop:#1}}
\newcommand{\reco}[1]{Corollary~\protect\ref{cor:#1}}
\newcommand{\resec}[1]{Section~\protect\ref{sec:#1}}

\newcommand{\rerems}[1]{Remarks~\protect\ref{rem:#1}}
\newcommand{\req}[1]{equation~(\protect\ref{eq:#1})}
\newcommand{\reqref}[1]{(\protect\ref{eq:#1})}

\begin{document}

\title{The classification and the conjugacy classes of the finite subgroups of the sphere braid groups}

\author{DACIBERG~LIMA~GON\c{C}ALVES\\
Departamento de Matem\'atica - IME-USP,\\
Caixa Postal~66281~-~Ag.~Cidade de S\~ao Paulo,\\ 
CEP:~05311-970 - S\~ao Paulo - SP - Brazil.\\
e-mail:~\texttt{dlgoncal@ime.usp.br}\vspace*{4mm}\\
JOHN~GUASCHI\\
Laboratoire de Math\'ematiques Nicolas Oresme UMR CNRS~\textup{6139},\\
Universit\'e de Caen BP 5186,\\
14032 Caen Cedex, France.\\
e-mail:~\texttt{guaschi@math.unicaen.fr}}

\date{21st November 2007}

\begingroup
\renewcommand{\thefootnote}{}
\footnotetext{
(secondary).}
\endgroup 


\maketitle

\begin{abstract}\noindent
\emph{Let $n\geq 3$. We classify the finite groups which are realised as subgroups of the sphere
braid group $B_n(\St)$. Such groups must be of cohomological period $2$ or $4$. Depending on the
value of $n$, we show that the following are the maximal finite subgroups of $B_n(\St)$:
$\Z_{2(n-1)}$; the dicyclic groups of order $4n$ and $4(n-2)$; the binary tetrahedral group $T_1$;
the binary octahedral group $O_1$; and the binary icosahedral group $I$. We give geometric as well
as some explicit algebraic constructions of these groups in $B_{n}(\St)$, and determine the number
of conjugacy classes of such finite subgroups. We also reprove Murasugi's classification of the
torsion elements of $B_{n}(\St)$, and explain how the finite subgroups of $B_{n}(\St)$ are related
to this classification, as well as to the lower central and derived series of $B_{n}(\St)$.}
\end{abstract}

\section{Introduction}\label{sec:intro}
 
The braid groups $B_n$ of the plane were introduced by E.~Artin in~1925~\cite{A1,A2}. Braid groups
of surfaces were studied by Zariski~\cite{Z}. They were later generalised 
by Fox to braid groups of arbitrary topological spaces via the following definition~\cite{FoN}. Let
$M$ be a compact, connected surface, and let $n\in\N$. We denote the set of all ordered $n$-tuples
of distinct points of $M$, known as the \emph{$n\th$ configuration space of $M$}, by:
\begin{equation*}
F_n(M)=\setr{(p_1,\ldots,p_n)}{\text{$p_i\in M$ and $p_i\neq p_j$ if $i\neq j$}}.
\end{equation*}
Configuration spaces play an important r\^ole in several branches of mathematics and have been
extensively studied, see~\cite{CG,FH} for example. 

The symmetric group $\sn$ on $n$ letters acts freely on $F_n(M)$ by
permuting coordinates. The corresponding quotient will be denoted by
$D_n(M)$. The \emph{$n\th$ pure braid group $P_n(M)$} (respectively
the \emph{$n\th$ braid group $B_n(M)$}) is defined to be the
fundamental group of $F_n(M)$ (respectively of $D_n(M)$). 

Together with the real projective plane $\rp$, the braid groups of the
$2$-sphere $\St$ are of particular interest, notably because they have
non-trivial centre~\cite{GVB,GG2}, and torsion elements~\cite{vB,M}.
Indeed, Van Buskirk showed that among the braid groups of compact,
connected surfaces, $B_n(\St)$ and $B_n(\rp)$ are the only ones to
have torsion~\cite{vB}. Let us recall briefly some of the properties
of $B_n(\St)$~\cite{FvB,GVB,vB}.

If $\dt\subseteq \St$ is a topological disc, there is a group
homomorphism $\map {\iota}{B_n}[B_n(\St)]$ induced by the inclusion.
If $\beta\in B_n$, we shall denote its image $\iota(\beta)$ simply
by $\beta$. Then $B_n(\St)$ is generated by
$\sigma_1,\ldots,\sigma_{n-1}$ which are subject to the following
relations:
\begin{equation*}\label{eq:presnbns}
\begin{gathered}
\text{$\sigma_{i}\sigma_{j}=\sigma_{j}\sigma_{i}$ if $\lvert i-j\rvert\geq 2$ and $1\leq i,j\leq n-1$}\\
\text{$\sigma_{i}\sigma_{i+1}\sigma_{i}=\sigma_{i+1}\sigma_{i}\sigma_{i+1}$ for all $1\leq i\leq n-2$, and}\quad\\
\text{$\sigma_1\cdots \sigma_{n-2}\sigma_{n-1}^2 \sigma_{n-2}\cdots \sigma_1=1$.}
\end{gathered}
\end{equation*}
Consequently, $B_n(\St)$ is a quotient of $B_n$. The first three sphere braid groups are finite: $B_1(\St)$ is trivial,
$B_2(\St)$ is cyclic of order~$2$, and $B_3(\St)$ is a $\text{ZS}$-metacyclic group (a group whose Sylow subgroups,
commutator subgroup and commutator quotient group are all cyclic) of order~$12$, isomorphic to the semi-direct product
$\Z_3 \rtimes \Z_4$ of cyclic groups, the action being the non-trivial one, which in turn is isomorphic to the dicyclic
group $\dic{12}$ of order~$12$. The Abelianisation of $B_n(\St)$ is isomorphic to the cyclic group $\Z_{2(n-1)}$. The
kernel of the associated projection $\map{\xi}{B_n(\St)}[\Z_{2(n-1)}]$ (which is defined by $\xi(\sigma_i)=
\overline{1}$ for all $1\leq i\leq n-1$) is the commutator subgroup $\Gamma_2\left(B_n(\St) \right)$. If $w\in B_n(\St)$
then $\xi(w)$ is  the exponent sum (relative to the $\sigma_i$) of $w$ modulo $2(n-1)$. 

Gillette and Van Buskirk showed that if $n\geq 3$ and $k\in \N$ then
$B_n(\St)$ has an element of order $k$ if and only if $k$ divides one
of $2n$, $2(n-1)$ or $2(n-2)$~\cite{GVB}. The torsion elements of
$B_n(\St)$ and $B_n(\rp)$ were later characterised by
Murasugi~\cite{M}. For $B_n(\St)$, these elements are as follows: 
\begin{thm}[\cite{M}]\label{th:murasugi}
Let $n\geq 3$. Then the torsion elements of $B_n(\St)$ are precisely powers of conjugates of the following three elements:
\begin{enumerate}[(a)]
\item $\alpha_0= \sigma_1\cdots \sigma_{n-2}
\sigma_{n-1}$ (which is of order $2n$).
\item $\alpha_1=\sigma_1\cdots
\sigma_{n-2} \sigma_{n-1}^2$ (of order $2(n-1)$).
\item $\alpha_2=\sigma_1\cdots \sigma_{n-3} \sigma_{n-2}^2$ (of order
$2(n-2)$).
\end{enumerate}
\end{thm}
The three elements $\alpha_0$, $\alpha_1$ and $\alpha_2$ are
respectively $n\th$, $(n-1)\th$ and $(n-2)\th$ roots of $\ft$, where
$\ft$ is the so-called `full twist' braid of $B_n(\St)$, defined by
$\ft= (\sigma_1\cdots\sigma_{n-1})^n$. So $B_n(\St)$ admits finite
cyclic subgroups isomorphic to $\Z_{2n}$, $\Z_{2(n-1)}$ and
$\Z_{2(n-2)}$. In~\cite{GG4}, we showed that $B_n(\St)$ is generated
by $\alpha_0$ and $\alpha_1$. If $n\geq 3$, $\ft[n]$ is the unique
element of $B_n(\St)$ of order $2$, and it generates the centre of
$B_n(\St)$. It is also the square of the \emph{Garside element} (or
`half twist') defined by:
\begin{equation*}
\garside = (\sigma_1 \cdots \sigma_{n-1}) (\sigma_1 \cdots \sigma_{n-2}) \cdots (\sigma_1 \sigma_2)\sigma_1.
\end{equation*}

For $n\geq 4$, $B_n(\St)$ is infinite. It is an interesting question as to which
finite groups are realised as subgroups of $B_n(\St)$ (apart of course from the
cyclic groups $\ang{\alpha_i}$ and their subgroups given in \reth{murasugi}). Another question is the
following: how many conjugacy classes are there in $B_n(\St)$ of a given
abstract finite group? As a partial answer to the first question, we proved
in~\cite{GG4} that $B_n(\St)$ contains an isomorphic copy of the finite group
$B_3(\St)$ of order $12$ if and only if $n\not\equiv 1 \bmod 3$. 

While studying the lower central and derived series of the sphere
braid groups, we showed that $\Gamma_2\left(B_4(\St) \right)$ is
isomorphic to a semi-direct product of $\quat$ by a free group of rank
$2$~\cite{GG5}. After having proved this result, we noticed that the
question of the realisation of $\quat$ as a subgroup of $B_n(\St)$ had been
explicitly posed by R.~Brown~\cite{DD} in connection with the Dirac
string trick~\cite{Fa,Ne} and the fact that the fundamental group of
$\operatorname{SO}(3)$ is isomorphic to $\Z_2$. The case $n=4$ was
studied by J.~G.~Thompson~\cite{Thp}. In a previous paper, we provided a complete answer to this question:
\begin{thm}[\cite{GG6}]\label{th:quat8}
Let $n\in\N$, $n\geq 3$.
\begin{enumerate}
\item $B_n(\St)$ contains a subgroup isomorphic to $\quat$ if and only if $n$ is even.
\item If $n$ is divisible by $4$ then $\Gamma_2\left(B_n(\St) \right)$ contains a subgroup isomorphic to $\quat$.
\end{enumerate}
\end{thm}
As we also pointed out in~\cite{GG6}, for all $n\geq 3$, the construction of $\quat$ may be generalised in order to
obtain a subgroup $\ang{\alpha_0, \garside}$ of $B_n(\St)$ isomorphic to the dicyclic group $\dic{4n}$ of order $4n$.

It is thus natural to ask which other finite groups are realised as
subgroups of $B_n(\St)$. One common property of the above subgroups is
that they are finite periodic groups of cohomological period $2$ or
$4$. In fact, this is true for all finite subgroups of $B_n(\St)$.
Indeed, by~\cite{GG4}, the universal covering $X$ of $F_n(\St)$ is a
finite-dimensional complex which has the homotopy type of $\St[3]$ (we were recently informed by V.~Lin that $X$ is
biholomorphic to the direct product of $\operatorname{SL}(2,\mathbb{C})$ by the Teichm\"uller space of the $n$-punctured
Riemann sphere~\cite{Li}). Thus any finite subgroup of $B_n(\St)$ acts freely on $X$, and so has
period $2$ or $4$ by Proposition~10.2, Section~10, Chapter~VII
of~\cite{Br}. Since $\ft$ is the unique element of order~$2$ of
$B_n(\St)$, and it generates the centre $Z(B_n(\St))$, the Milnor
property must be satisfied for any finite subgroup of $B_n(\St)$. 
Recall also that a finite periodic group $G$ satisfies the $p^2$-condition (if $p$ is prime and divides the order of $G$ then $G$ has
no subgroup isomorphic to $\Z_p \times \Z_p$), which implies that a
Sylow $p$-subgroup of $G$ is cyclic or generalised quaternion, as well
as the $2p$-condition (each subgroup of order $2p$ is cyclic). The
classification of finite periodic groups is given by the
Suzuki-Zassenhaus theorem (see~\cite{AM,Th} for example), and thus
provides a possible line of attack for the subgroup realisation
problem. The periods of the different families of these groups were
determined in a series of papers by Golasi\'nski and
Gonçalves~\cite{gg1,gg2,gg3,gg4,gg5, gg6}, and so in theory we may obtain a
list of those of period~$4$. A list of all periodic groups of period~$4$
is provided in~\cite{Th}. However, in the current context, a more
direct approach is obtained via the relationship between the braid
groups and the mapping class groups of $\St$, which we shall now
recall.

For $n\in\N$, let $\mathcal{M}_{0,n}$ denote the mapping class group of the $n$-punctured sphere. We allow the $n$
marked points to be permuted. If $n\geq 2$, a presentation of $\mathcal{M}_{0,n}$ is obtained from that of $B_n(\St)$
by adding the relation $\ft=1$~\cite{Ma,MKS}. In other words, we have the following central extension:
\begin{equation}\label{eq:mcg}
1 \to \ang{\ft} \to B_n(\St) \stackrel{p}{\to} \mathcal{M}_{0,n} \to 1.
\end{equation}
If $n=2$, $B_2(\St) \cong \mathcal{M}_{0,2}\cong \Z_2$. For $n=3$, since $\mathcal{M}_{0,3}\cong \sn[3]$, this short
exact sequence does not split, and in fact for $n\geq 4$ it does not split either~\cite{GVB}.

This exact sequence may also be obtained in the following
manner~\cite{Bi}. Let $\diff{\St}$ denote the group of
orientation-preserving homeomorphisms of $\St$, and let $X \in
D_n(\St)$. Then $\diff{\St,X}=\setl{f\in \diff{\St}}{f(X)=X}$ is a
subgroup of $\diff{\St}$, and we have a fibration $\diff{\St,X}\to
\diff{\St}\to D_n(\St)$, where the basepoint of $D_n(\St)$ is taken to
be $X$, and where the second map evaluates an element of $\diff{\St}$
on $X$. The resulting long exact sequence in homotopy yields:
\begin{multline}\label{eq:longexact}
\cdots \to \pi_1\left(\diff{\St,X} \right) \to \underbrace{\pi_1\left(\diff{\St} \right)}_{\Z_2}  \to \underbrace{\pi_1(D_n(\St))}_{B_n(\St)}\\ 
 \stackrel{\delta}{\to} \underbrace{\pi_0\left(\diff{\St,X} \right)}_{\mathcal{M}_{0,n}} \to \underbrace{\pi_0\left(\diff{\St} \right)}_{=\brak{1}}.
\end{multline}
The homomorphism $\map{\delta}{B_n(\St)}[\mathcal{M}_{0,n}]$ is the boundary
operator which we shall use in \resec{geomreal} in order to describe the
geometric realisation of the finite subgroups of $B_n(\St)$. If $n\geq 3$ then
$\pi_1\left(\diff{\St,X} \right)=\brak{1}$~\cite{EE,Hm,Sc}, and we thus
recover~\req{mcg} (the interpretation of the Dirac string trick in terms of the
sphere braid groups~\cite{Fa,H,Ne} gives rise to the identification of
$\pi_1\left(\diff{\St} \right)$ with $\ang{\ft}$).

In a recent paper, Stukow applies Kerckhoff's solution of the Nielsen realisation problem~\cite{K} to classify the
finite maximal subgroups of $\mathcal{M}_{0,n}$~\cite{S}. Applying his results to~\req{mcg}, we shall see in
\resec{class} that their counterparts in $B_n(\St)$ are cyclic, dicyclic and binary polyhedral groups:
\begin{thm}\label{th:finitebn}
Let $n\geq 3$. The maximal finite subgroups of $B_n(\St)$ are:
\begin{enumerate}[(a)]
\item\label{it:fina} $\Z_{2(n-1)}$ if $n\geq 5$.
\item\label{it:finb} the dicyclic group $\dic{4n}$ of order $4n$.
\item\label{it:finc} the dicyclic group $\dic{4(n-2)}$ if $n=5$ or $n\geq 7$.
\item\label{it:find} the binary tetrahedral group, denoted by $T_1$, if $n\equiv 4 \bmod 6$.
\item\label{it:fine} the binary octahedral group, denoted by $O_1$, if $n\equiv 0,2\bmod 6$.
\item\label{it:finf} the binary icosahedral group, denoted by $I$, if $n\equiv 0,2,12,20\bmod 30$.
\end{enumerate}
\end{thm}

\begin{rems}\mbox{}
\begin{enumerate}[(a)]
\item If $n$ is odd then the only finite subgroups of $B_n(\St)$ are cyclic or dicyclic. In the latter case, the
dicyclic group $\dic{4n}$ (resp.\ $\dic{4(n-2)}$) is ZS-metacyclic~\cite{CM}, and is isomorphic to $\Z_n\rtimes \Z_4$
(resp.\ $\Z_{n-2} \rtimes \Z_4$), where the action is multiplication by~$-1$.
\item If $n$ is even then one of the binary tetrahedral or octahedral groups is realised as a maximal finite subgroup of
$B_n(\St)$. Further, since $T_1$ is a subgroup of $O_1$, $T_1$ is realised as a subgroup of $B_n(\St)$ for all $n$ even,
$n\geq 4$.
\item The groups of \reth{finitebn} and their subgroups are the finite groups of quaternions~\cite{Co}. Indeed, for $p,q,r \in \N$, let us denote
\begin{equation*}
\ang{p,q,r}=\setangr{A,B,C}{A^p=B^q= C^r=ABC}.
\end{equation*}
Then $\Z_{2(n-1)}= \ang{n-1,n-1,1}$, $\dic{4n}= \ang{n,2,2}$, $\dic{4(n-2)}= \ang{n-2,2,2}$, $T_1=\ang{3,3,2}$, $O_1=
\ang{4,3,2}$ and $I= \ang{5,3,2}$. It is shown in~\cite{Co,CM} that for $T_1$, $O_1$ and $I$, this presentation is
equivalent to:
\begin{equation*}
\ang{p,3,2}=\setangr{A,B}{A^p=B^3= (AB)^2},
\end{equation*}
for $p\in \brak{3,4,5}$, and that the element $A^p$ is central and is the unique element of order~$2$ of
$\ang{p,3,2}$.
\end{enumerate}
\end{rems}

In \resec{class}, we also generalise another result of Stukow concerning the conjugacy classes of finite subgroups of
$\mathcal{M}_{0,n}$ to $B_n(\St)$:
\begin{prop}\label{prop:isoconj}\mbox{}
\begin{enumerate}[(a)]
\item\label{it:maxconj} Two maximal finite subgroups of $B_n(\St)$ are isomorphic if and only if they are conjugate.
\item\label{it:conj2} Each abstract finite subgroup $G$ of $B_n(\St)$ is realised as a single conjugacy class within
$B_n(\St)$, with the exception, when $n$ is even, of the following cases, for which there are precisely two conjugacy
classes:
\begin{enumerate}[(i)]
\item\label{it:z4} $G=\Z_4$.
\item\label{it:dic2m} $G=\dic{4r}$, where $r$ divides $\frac{n}{2}$ or $\frac{n-2}{2}$.
\end{enumerate}
\end{enumerate}
\end{prop}

In \resec{geomreal}, we explain how to obtain geometrically the subgroups of \reth{finitebn}, and we also give explicit
group presentations of the cyclic and dicyclic subgroups, as well as in the special case $T_{1}$ for $n=4$. 

In order to understand better the finite subgroups of $B_{n}(\St)$, it is often useful to know their relationship with
the three classes of elements described in \reth{murasugi}. This shall be carried out in \repr{classgi} (see \resec{muragi}).

The two conjugacy classes of part~(\ref{it:conj2})(\ref{it:z4}) are
realised by the subgroups $\ang{\alpha_0^{n/2}}$ and
$\ang{\alpha_2^{(n-2)/2}}$ (they are non conjugate since they project
to non-conjugate subgroups in $\sn$). In \resec{gamma2}, we
construct the two conjugacy classes of
part~(\ref{it:conj2})(\ref{it:dic2m}) of \repr{isoconj}:
\begin{thm}\label{th:gammaconj}
Let $n\geq 4$ be even. Let $N\in \brak{n,n-2}$, and let $x=\alpha_0$ (resp.\ $x=\alpha_0 \alpha_2 \alpha_0^{-1}$) if
$N=n$ (resp.\ $N=n-2$). Set $N=2^l k$, where $l\in \N$, and $k$ is odd. Then for $j=0,1,\ldots, l$, and $q$ a divisor of
$k$, we have:
\begin{enumerate}[(a)]
\item\label{it:2jcopies} $B_n(\St)$ contains $2^j$ copies of $\dic{2^{l+2-j}k/q}$ of the form $\ang{x^{2^jq}, x^{iq}\garside}$, where $i=0,1, \ldots, 2^j-1$.
\item\label{it:conjq8} if $0\leq i,i'\leq 2^j-1$, $\ang{x^{2^jq}, x^{iq}\garside}$ and $\ang{x^{2^jq}, x^{i'q}\garside}$ are conjugate if and only if $i-i'$ is even.
\end{enumerate}
\end{thm}

Another question arising from \reth{quat8} is the existence of copies of $\quat$ lying in $\Gamma_2(B_n(\St))$. More
generally, one may ask whether the dicyclic groups constructed above (and indeed the other finite subgroups of
$B_{n}(\St)$) are contained in $\Gamma_2(B_n(\St))$. In the dicyclic case, we have the following result, also proved
in \resec{gamma2}:
\begin{prop}\label{prop:quatgamma2}
Let $n\geq 4$ be even, let $N\in \brak{n,n-2}$, and let $r$ divide $N$. If $r$ does not divide $N/2$ then the
subgroups of $B_{n}(\St)$ abstractly isomorphic to $\dic{4r}$  are not contained in $\Gamma_{2}(B_{n}(\St))$. If $r$
divides $N/2$ then up to conjugacy, $B_{n}(\St)$ has a two subgroups abstractly isomorphic to $\dic{4r}$, one of which
is contained in $\Gamma_{2}(B_{n}(\St))$, and the other not. In particular, $B_n(\St)$ exhibits the two conjugacy
classes of $\quat$, one of which lies in $\Gamma_2(B_n(\St))$, the other not. 
\end{prop}
The corresponding result for the binary polyhedral groups may be found in \repr{gamma2}. As a
corollary of our results we obtain an alternative proof of \reth{murasugi} (see \resec{murasugi}).

\subsection*{Acknowledgements}

This work took place during the visits of the second author to the Departmento de Matem\'atica do IME-Universidade de
S\~ao Paulo during the periods 10\up{th}~--~27\up{th}~August 2007 and 4\up{th}~--~23\up{rd}~October 2007,  and of the
visit of the first author to the Laboratoire de Math\'ematiques Emile Picard, Université Paul Sabatier, Toulouse, during
the period 27\up{th}~August~--~17\up{th}~September 2007. It was supported by the international Cooperation USP/Cofecub
project number 105/06, by an `Aide Ponctuelle de Coop\'eration' from the Universit\'e Paul Sabatier, by the
Pr\'o-Reitoria de Pesquisa Projeto~1, and by the Projeto Tem\'atico FAPESP no.~2004/10229-6 `Topologia Alg\'ebrica,
Geometrica e Diferencial'.

\section{The classification of the finite maximal subgroups of $B_n(\St)$}\label{sec:class}

In this section, we prove \reth{finitebn}. We start by making some remarks concerning the central
extension~(\ref{eq:mcg}).

\begin{rems}\label{rem:bnmcg}
Let $G$ be a finite subgroup of $B_n(\St)$.
\begin{enumerate}[(a)]
\item If $H$ is a finite subgroup of $\mathcal{M}_{0,n}$ then $p^{-1}(H)$ is a finite subgroup of $B_n(\St)$ of order $2\,\lvert H\rvert$.
\item\label{it:odd} If $\lvert G\rvert$ is odd then $\ft \notin G$, and so $G \cong p(G)$. Conversely, if $G \cong p(G)$
then $p\vert_G$ is injective, and thus $\ft \notin G$, so $\lvert G\rvert$ is odd.
\item\label{it:even} If $\lvert G\rvert$ is even then $\ft \in G$, and so we obtain the following short exact sequence:
\begin{equation}\label{eq:mcgfinite}
1 \to \ang{\ft} \to G \stackrel{p\vert_{G}}{\to} p(G) \to 1,
\end{equation}
where $p(G)$ is a finite subgroup of $\mathcal{M}_{0,n}$ of order $\displaystyle \frac{\lvert G\rvert}{2}$. 
\item If $G$ is a maximal finite subgroup of $B_n(\St)$ then $\lvert G\rvert$ is even, and $p(G)$ is a maximal finite
subgroup of $\mathcal{M}_{0,n}$. Conversely, if $H$ is a maximal finite subgroup of $\mathcal{M}_{0,n}$ then $p^{-1}(H)$
is a maximal finite subgroup of $B_n(\St)$.
\end{enumerate}
\end{rems}

We recall Stukow's theorem:
\begin{thm}[\cite{S}]\label{th:stukow}
Let $n\geq 3$. The maximal finite subgroups of $\mathcal{M}_{0,n}$ are:
\begin{enumerate}[(a)]
\item\label{it:stuka} $\Z_{n-1}$ if $n\neq 4$.
\item the dihedral group $\dih{2n}$ of order $2n$.
\item the dihedral group $\dih{2(n-2)}$ if $n=5$ or $n\geq 7$.
\item $\an[4]$ if $n\equiv 4,10 \bmod 12$.
\item $\sn[4]$ if $n\equiv 0,2,6,8,12,14,18,20 \bmod 24$.
\item $\an[5]$ if $n\equiv 0,2,12,20,30,32,42,50\bmod 60$.
\end{enumerate}
\end{thm}

\begin{rem}
In the case $n=3$, $\mathcal{M}_{0,3}$ is isomorphic to $\dih{6}$, obtained as a maximal subgroup in part~(b) of
\reth{stukow}, and so its subgroup isomorphic to $\Z_{2}$ is not maximal. This explains the discrepancy between the
value of $n$ in part~(a) of Theorems~\ref{th:finitebn} and~\ref{th:stukow}.
\end{rem}

\begin{proof}[Proof of \reth{finitebn}]
By \rerems{bnmcg}, we just need to check that the given groups are those obtained as extensions of $\ang{\ft}$ by the
groups of \reth{stukow}. We start by making some preliminary remarks. Let $H$ be one of the finite maximal subgroups of
$\mathcal{M}_{0,n}$, and let $G$ be a finite (maximal) subgroup of $B_n(\St)$ of order $2\, \lvert H \rvert$ which fits
into the following short exact sequence:
\begin{equation}\label{eq:basic}
1 \to \ang{\ft} \to G \stackrel{p\vert_{G}}{\to} H \to 1,
\end{equation}
where $\ft\in G$ belongs to the centre of $G$, and is the unique element of $G$ of order~$2$. Then $G=p^{-1}(H)$, and so is unique.

Suppose that $y\in H$ is of order $k\geq 2$. Then $y$ has two preimages in $G$, of the form $x$ and $x\ft$, say, and $x$
is of order $k$ or $2k$. If $k$ is even then by \rerems{bnmcg}(\ref{it:even}), $x$ must be of order $2k$, $x^k=\ft$ and
$\ft \in \ang{x}$. If $k$ is odd then $x$ is of order $k$ (resp.\ $2k$) if and only if $x\ft$ is of order $2k$ (resp.\
$k$).

A presentation of $G$ may be obtained by applying standard results
concerning the presentation of an extension (see Theorem~1, Chapter~13
of~\cite{J}). If $H$ is generated by $h_1,\ldots, h_k$ then $G$ is
generated by $g_1,\ldots, g_k,\ft$, where $p(g_i)=h_i$ for
$i=1,\ldots,k$. One relation of $G$ is just $\ft^2=1$, that of
$\ker{p}$. Since $\ker{p}\subseteq Z(G)$, the remaining relations of
$G$ are obtained by rewriting the relators of $H$ in terms of the
coset representatives, and expressing the corresponding element in the
form $\ft^{\epsilon}$, where $\epsilon\in \brak{0,1}$.

We consider the six cases of \reth{stukow} as follows.
\begin{enumerate}[(a)]
\item $H \cong \Z_{n-1}$: let $y$ be a generator of $H$, and let $x\in G$ be such that $p(x)=y$. Then $G=\ang{\ft,x}$
and $\lvert G \rvert= 2(n-1)$. If $n$ is odd then $\ft\in \ang{x}$, $G=\ang{x}$,  and $x$ is of order $2(n-1)$. If $n$
is even then $G=\ang{x\ft}$ (resp.\ $G=\ang{x}$) if $x$ is of order $n-1$ (resp.\ $2(n-1)$), and $ G \cong \Z_{2(n-1)}$
in both cases. 
\item $H \cong \dih{2n}$: let $y,z\in H$ be such that $o(y)=n$, $o(z)=2$ and $zyz^{-1}=y^{-1}$, and let $x,w\in G$ be
such that $p(x)=y$ and $p(w)=z$. So $G=\ang{\ft,x,w}$ and $\lvert G \rvert=4n$. 
From above, it follows that $w^2=\ft$, so $G=\ang{x,w}$. If $n$ is even then $x$ is of order $2n$ and $x^n=\ft$. The
same result may be obtained if $n$ is odd, replacing $x$ by $x\ft$ if necessary. Further, $wxw^{-1}x \in\ker{p}$. If
$wxw^{-1}x=\ft$ then $(wx)^2=1$. So either $w=x^{-1}$ or $wx=\ft$, and in both cases we conclude that $G=\ang{x}$ which
contradicts $\lvert G \rvert=4n$. Hence $wxw^{-1}x=1$, and since $\lvert G \rvert=4n$, $G$ is isomorphic to $\dic{4n}$.
\item $H \cong \dih{2(n-2)}$: the previous argument shows that $G\cong \dic{4(n-2)}$.

\item Suppose that $H$ is isomorphic to one of the remaining groups $\an[4]$, $\sn[4]$ or $\an[5]$ of \reth{stukow}. Let
$p=3$ if $H\cong \an[4]$, $p=4$ if $H\cong \sn[4]$, and $p=5$ if $H\cong \an[5]$. Then $H$ has a presentation given
by~\cite{Co,CM}:
\begin{equation*}
H=\setangr{u,v}{u^2=v^3=(uv)^p=1}.
\end{equation*}
Let $x,w\in G$ be such that $p(x)=u$ and $p(w)=v$. Then $G=\ang{x,w,\ft}$. From above, we must have $x^2=\ft$. Further,
replacing $w$ by $w\ft$, we may suppose that $w^3=\ft$. If $p=4$ then $(xw)^p=\ft$, while if $p\in \brak{3,5}$,
replacing $x$ by $x\ft$ if necessary, we may suppose that $(xw)^p=\ft$. It is shown in~\cite{Co,CM} that
$x^2=w^3=(xw)^p=\ft$ implies that $\ft^2=1$, so $G$ admits a presentation given by:
\begin{equation*}
G=\setangr{x,w}{x^2=w^3=(xw)^p}.
\end{equation*}
Thus $G\cong T_1$ if $p=3$, $G\cong O_1$ if $p=4$ and $G\cong I$ if $p=5$. This completes the proof of the theorem.
\end{enumerate}
\end{proof}

\begin{rems}\label{rem:isoconj}
Let $G_1,G_2$ be finite subgroups of $B_n(\St)$. 
\begin{enumerate}[(a)]
\item If they are of odd order then by \rerems{bnmcg}, $G_1$ and $G_2$
are isomorphic if and only if $p(G_1)$ and $p(G_2)$ are isomorphic. So
suppose that $G_1$ and $G_2$ are of even order. If $p(G_1)$ and
$p(G_2)$ are isomorphic then it follows from the construction of
\reth{finitebn} that $G_1$ and $G_2$ are isomorphic. Conversely,
suppose that $G_1$ and $G_2$ are isomorphic via an isomorphism
$\map{\alpha}{G_1}[G_2]$. Since $\ft$ belongs to both, and is the
unique element of order~$2$, we must have $\alpha(\ft)=\ft$, and thus
$\alpha$ induces an isomorphism
$\map{\widetilde{\alpha}}{p(G_1)}[p(G_2)]$ satisfying
$\widetilde{\alpha}\circ p =p\circ \alpha$.
\item\label{it:conj} If $G_1,G_2$ are conjugate then clearly so are $p(G_1)$ and
$p(G_2)$. Conversely, suppose that $p(G_1), p(G_2)$ are conjugate subgroups of
$\mathcal{M}_{0,n}$. Then there exists $g\in \mathcal{M}_{0,n}$ such that
$p(G_2)=g p(G_1) g^{-1}$. If $G_1$ and $G_2$ are of even order, the fact
that~\req{mcg} is a central extension implies that $G_1,G_2$ are conjugate.
If $G_1$ and $G_2$ are of odd order, let $L_i=p^{-1}(p(G_i))$ for $i=1,2$. Then
$[L_i:G_i]=2$, and it follows from the even order case that $L_1$ and $L_2$ are
conjugate in $B_n(\St)$. But $L_i=G_i \coprod \ft G_i$, and its odd order
elements are precisely those of $G_i$. So the conjugacy between $L_1$ and $L_2$
must send $G_1$ onto $G_2$.
\end{enumerate}
\end{rems}

We are now able to prove \repr{isoconj}.

\begin{proof}[Proof of \repr{isoconj}]
Part~(\ref{it:maxconj}) follows from \rerems{bnmcg} and~\ref{rem:isoconj}. To prove part~(\ref{it:conj2}), let $G_1,G_2$
be abstractly isomorphic finite subgroups of $B_n(\St)$, and for $i=1,2$, let $H_i=p(G_i)$. Then $H_1\cong H_2$: if the
$G_i$ are of odd order then $H_i\cong G_i$, so $H_1 \cong H_2$, while if the $G_i$ are of even order, any isomorphism
between them must send $\ft \in G_1$ onto $\ft \in G_2$, and so projects to an isomorphism between the $H_i$.
From \rerems{isoconj}(\ref{it:conj}), $G_1$ and $G_2$ are conjugate if and only if $H_1$ and $H_2$ are, and so the
number of conjugacy classes of subgroups of $B_n(\St)$ isomorphic to $G_1$ is the same as the number of conjugacy
classes of subgroups of $\mathcal{M}_{0,n}$ isomorphic to $H_1$. The result follows from the proof of \reth{finitebn} by
remarking that a subgroup of $\mathcal{M}_{0,n}$ isomorphic to $\Z_2$ (resp.\ $\dih{2r}$) lifts to a subgroup of
$B_n(\St)$ which is isomorphic to $\Z_4$ (resp.\ $\dic{4r}$).
\end{proof}

\section{Realisation of the maximal finite subgroups of $B_n(\St)$}\label{sec:geomreal}

In this section, we analyse the geometric and algebraic realisations of the subgroups given in \reth{finitebn}.

\subsection{The algebraic realisation of some finite subgroups of $B_{n}(\St)$}

The maximal cyclic and dicyclic subgroups of $B_n(\St)$ may be realised as follows:
\begin{enumerate}[(a)]
\item $\Z_{2(n-1)}\cong \ang{\alpha_1}$.
\item $\dic{4n}\cong \ang{\alpha_0, \garside}$~\cite{GG6}. 
\item The algebraic realisation of $\dic{4(n-2)}$ is given by the following proposition:

\begin{prop}
For all $n\geq 3$, the subgroup $\ang{\alpha_0 \alpha_2 \alpha_0^{-1}, \garside}$ of $B_n(\St)$ is isomorphic to $\dic{4(n-2)}$.
\end{prop}

\begin{proof}
Let $x=\alpha_0 \alpha_2 \alpha_0^{-1}$. We know that $x$ is of order $2(n-2)$, and that $x^{n-1}=\ft=\garside^2$.
Further, by standard properties of the corresponding elements in $B_n$~\cite{Bi}, $\alpha_0 \sigma_i \alpha_0^{-1}=
\sigma_{i+1}$ for all $i=1,\ldots, n-2$, and $\garside \sigma_i \garside^{-1}= \sigma_{n-i}$ for all $i=1,\ldots, n-1$.
Hence $x= \sigma_2 \cdots \sigma_{n-2} \sigma_{n-1}^2$, and
\begin{equation*}
\garside x \garside^{-1} = \sigma_{n-2} \cdots \sigma_2 \sigma_1^2 = \sigma_{n-1}^{-2} \sigma_{n-2}^{-1} \cdots \sigma_2^{-1}= x^{-1}.
\end{equation*}
Thus $\ang{x, \garside}$ is isomorphic to a quotient of $\dic{4(n-2)}$. But $\garside \notin \ang{x}$, so $\ang{x,
\garside}$ contains the $2(n-2)+1$ distinct elements of $\ang{x}\cup \brak{\garside}$, and the result follows.
\end{proof}

\begin{rem}
In the special case $n=4$, the binary tetrahedral group $T_1$ may be realised as follows. Let $y=\sigma_1\sigma_3^{-1}$.
From~\cite{GG6}, we know that $\ang{y, \garside[4]} \cong \quat$. In $B_4(\St)$, we also have $(\sigma_2\sigma_1)^3=
(\sigma_2\sigma_3)^3=\ft[4]=\garside[4]^2$. Then $\ang{\alpha_1^2}\cong \Z_3$ acts on $\ang{y, \garside[4]}$ as follows:
\begin{align*}
\alpha_1^2\cdot \garside[4] \cdot \alpha_1^{-2} &= \alpha_1^2 (\garside[4] \alpha_1^{-2}\garside[4]^{-1}) \garside[4]\\
&= \alpha_1^2 (\sigma_1^{-2}\sigma_2^{-1}\sigma_3^{-1})^2 \garside[4] \;\text{(by the action of $\garside[4]$)}\\
&= \alpha_1^2 (\sigma_2 \sigma_3)^2 \garside[4] \;\text{(using the surface relation of $B_{n}(\St)$)}\\
&= (\sigma_1\sigma_2 \sigma_3^2)^2  \cdot \sigma_3^{-1}\sigma_2^{-1}\cdot (\sigma_2 \sigma_3)^3 \garside[4]\\
&= \sigma_1\sigma_2 \sigma_3  \sigma_1\sigma_2 \sigma_1\cdot \sigma_1^{-1}\sigma_2^{-1} \sigma_1^{-1}\cdot \sigma_3
\sigma_1\sigma_2 \sigma_3\sigma_2^{-1} \garside[4]^3 \;\text{(as $\garside[4]^2=(\sigma_2\sigma_3)^3$)}\\
&= \garside[4] \sigma_1^{-1}\sigma_2^{-1}  \sigma_3 \sigma_2 \sigma_3\sigma_2^{-1} \garside[4]^3 \;\text{(as $\sigma_1$
commutes with $\sigma_3$)}\\
&=\garside[4] \sigma_1^{-1}\sigma_3 \garside[4]^3\;\text{(by the Artin braid relations)}\\
&=\garside[4]y^{-1}\garside[4]^{-1}=y \;\text{(by the action of $\garside[4]$ on $y$).}
\end{align*}
Further,
\begin{align*}
\alpha_1^2\cdot y \cdot \alpha_1^{-2} &= (\sigma_1^{-1}\sigma_2^{-1})^2 \cdot \sigma_1\sigma_3^{-1}\cdot (\sigma_2 \sigma_1)^2\\
&=(\sigma_1^{-1}\sigma_2^{-1})^2 \cdot \sigma_3^{-1}\sigma_2^{-1} \cdot (\sigma_2 \sigma_1)^3\;\text{(as $\sigma_1$ commutes
with $\sigma_3$)}\\
&= \sigma_1^{-1}\sigma_2^{-1}\sigma_1^{-1}\cdot \sigma_2^{-1} \sigma_3^{-1}\sigma_2^{-1}\cdot \garside[4]^2 \;\text{(as $\garside[4]^2=(\sigma_2\sigma_1)^3$)}\\
&= \sigma_1^{-1}\sigma_2^{-1}\sigma_1^{-1} \sigma_3^{-1}\sigma_2^{-1} \sigma_3^{-1}\cdot \garside[4]^2\;\text{(by the Artin braid relations)}\\
&= \sigma_1^{-1}\sigma_2^{-1}\sigma_1^{-1} \sigma_3^{-1}\sigma_2^{-1}\sigma_1^{-1}\cdot \sigma_1\sigma_3^{-1}\garside[4]^2\\
&= \garside[4]^{-1}y \garside[4]^2=\garside[4]y \;\text{(since $\garside[4]^2$ is central).}
\end{align*}
Hence $T_1=\quat \rtimes \Z_3\cong \ang{y,\garside[4]}\rtimes \ang{\alpha_1^2}$.
\end{rem}

\begin{rem}
We also have an algebraic representation of $T_{1}$ in $B_{6}(\St)$. Let
\begin{align*}
\gamma &= \sigma_{5} \sigma_{4} \sigma_{1}^{-1} \sigma_{2}^{-1}, \;\text{and} \\
\delta &= \sigma_{3}^{-1} \sigma_{4}^{-1} \sigma_{5}^{-1} \left ( \sigma_{2}^{-1} \sigma_{1}^{-1}
\sigma_{2}^{-1} \right) \sigma_{5} \sigma_{4} \sigma_{3}.
\end{align*}
Then we claim that $\ang{\gamma,\delta}\cong \quat\rtimes \Z_{3}\cong T_{1}$, where the action permutes the elements
$i,j,k$ of $\quat$. First, $\gamma^3=\delta^{2}=\ft[6]$. We now consider the subgroup $H=\ang{\delta,
\gamma \delta \gamma^{-1}}$. The action of conjugation by $\gamma$ permutes cyclically the elements $\delta$, $\gamma
\delta \gamma^{-1}$ and $\gamma \delta^2 \gamma^{-1}$, so is compatible with the action of $\Z_{3}$
on $\quat$. It just remains to show that $H\cong \quat$. Clearly $\delta^2=(\gamma \delta
\gamma^{-1})^2=\ft[6]$. Let us now prove that 
\begin{equation}\label{eq:relq8}
\delta^{-1} \cdot \gamma \delta \gamma^{-1} \cdot \delta= \gamma \delta^{-1} \gamma^{-1}. 
\end{equation}
Set $\rho=\sigma_{5} \sigma_{4} \sigma_{3}$, $\gamma'=\rho \gamma \rho^{-1}$ and $\delta'=\rho
\delta \rho^{-1}$. Then \req{relq8} is in turn equivalent to:
\begin{align*}
& \delta'^{-1} \cdot \gamma' \delta' \gamma'^{-1} \cdot \delta'= \gamma' \delta'^{-1}
\gamma'^{-1}\\
& \delta'^{-1} \gamma' \delta' \gamma'^{-1} \delta'^2 \delta'^{-1} \gamma' \delta' \gamma'^{-1}=1\\
& [\delta'^{-1}, \gamma']^2=\delta'^{-2}=\ft[6].
\end{align*}
We shall show that the latter relation holds. Notice that 
\begin{equation*}
\gamma'=\sigma_{5} \sigma_{4} \sigma_{3} \sigma_{5} \sigma_{4} \sigma_{1}^{-1} \sigma_{2}^{-1}
\sigma_{3}^{-1} \sigma_{4}^{-1} \sigma_{5}^{-1}= \sigma_{5} \sigma_{4} \sigma_{5} \sigma_{3} 
\sigma_{4} \alpha_{0}.
\end{equation*}
Then 
\begin{align*}
[\delta'^{-1}, \gamma']&= \sigma_{5}^{-1} \sigma_{4}^{-1} \sigma_{5}^{-1} \sigma_{2} \sigma_{1}
\sigma_{2} \cdot \sigma_{5} \sigma_{4} \sigma_{5} \sigma_{3} 
\sigma_{4} \alpha_{0} \cdot \sigma_{2}^{-1} \sigma_{1}^{-1} \sigma_{2}^{-1} \sigma_{5} \sigma_{4}
\sigma_{5} \cdot \alpha_{0}^{-1} \sigma_{4}^{-1} \sigma_{3}^{-1} \sigma_{5}^{-1} \sigma_{4}^{-1}
\sigma_{5}^{-1}\\
&= \sigma_{2} \alpha_{0} \sigma_{5}^{-1} \alpha_{0} \sigma_{2}^{-1} \sigma_{1}^{-1} \sigma_{2}^{-1}
\sigma_{5} \sigma_{4} \sigma_{5} \sigma_{5} \sigma_{4} \sigma_{3} \sigma_{2} \sigma_{1} \sigma_{4}^{-1} \sigma_{3}^{-1} \sigma_{5}^{-1} \sigma_{4}^{-1}
\sigma_{5}^{-1}\\
&= \sigma_{2} \alpha_{0} \sigma_{5}^{-1} \alpha_{0} \sigma_{2}^{-1} \sigma_{1}^{-1} \sigma_{2}^{-1}
\sigma_{5} \sigma_{3}^{-1} \sigma_{2}^{-1} \sigma_{1}^{-1} \sigma_{4}^{-1} \sigma_{3}^{-1}
\sigma_{5}^{-1} \sigma_{4}^{-1} \sigma_{5}^{-1}\\
&= \sigma_{2} \alpha_{0} \sigma_{5}^{-1} \alpha_{0} \sigma_{2}^{-1} \sigma_{5} \sigma_{1}^{-1}
\sigma_{2}^{-1} \sigma_{3}^{-1} \sigma_{4}^{-1} \sigma_{5}^{-1} \sigma_{1} \sigma_{1}^{-1} \sigma_{2}^{-1} \sigma_{1}^{-1} \sigma_{3}^{-1} \sigma_{4}^{-1} \sigma_{5}^{-1}\\
&= \sigma_{2} \alpha_{0} \sigma_{5}^{-1} \alpha_{0} \sigma_{2}^{-1} \sigma_{5} \alpha_{0} \sigma_{1} \sigma_{2}^{-1} \alpha_{0}\\
&= \sigma_{2} \alpha_{0} \sigma_{5}^{-1} \alpha_{0}^{-1} \alpha_{0}^{2} \sigma_{2}^{-1} \sigma_{5} \alpha_{0}^{-2}
\alpha_{0}^{3} \sigma_{1} \sigma_{2}^{-1} \alpha_{0}^{-3} \alpha_{0}^4\\
&= \sigma_{2} \tau^{-1} \sigma_{4}^{-1} \sigma_{1} \sigma_{4} \sigma_{5}^{-1} \alpha_{0}^4 = \sigma_{2} \tau^{-1}
\sigma_{1} \sigma_{5}^{-1} \alpha_{0}^4,
\end{align*}
since conjugation by $\alpha_{0}$ permutes cyclically the elements $\sigma_{1}, \sigma_{2}, \sigma_{3}, \sigma_{4},
\sigma_{5}$ and $\tau=\alpha_{0} \sigma_{5} \alpha_{0}^{-1}$. Thus
\begin{align*}
[\delta'^{-1}, \gamma']^2 &= \sigma_{2} \tau^{-1}
\sigma_{1} \sigma_{5}^{-1} \alpha_{0}^4 \sigma_{2} \tau^{-1}
\sigma_{1} \sigma_{5}^{-1} \alpha_{0}^{-4}\alpha_{0}^8= \sigma_{2} \tau^{-1}
\sigma_{1} \sigma_{5}^{-1} \tau \sigma_{4}^{-1} \sigma_{5} \sigma_{3}^{-1} \alpha_{0}^8.
\end{align*}
Let $\xi= \sigma_{2} \tau^{-1} \sigma_{1} \sigma_{5}^{-1} \tau \sigma_{4}^{-1} \sigma_{5} \sigma_{3}^{-1}$.
To prove that $[\delta'^{-1}, \gamma']^2= \ft[6]= \alpha_{0}^6$, it suffices to show that $\xi \alpha_{0}^2=1$. Now
\begin{align*}
\xi \alpha_{0}^2 &= \sigma_{2} \tau^{-1} \sigma_{1} \sigma_{5}^{-1} \tau \sigma_{4}^{-1} \sigma_{5} \sigma_{3}^{-1}
\alpha_{0}^2= \sigma_{2} \alpha_{0} \sigma_{5}^{-1} \alpha_{0}^{-1} \sigma_{1} \sigma_{5}^{-1} \alpha_{0} \sigma_{5} \alpha_{0}^{-1} \sigma_{4}^{-1} \sigma_{5}
\sigma_{3}^{-1} \alpha_{0}^2\\
&= \sigma_{2} \alpha_{0} \sigma_{5}^{-1} \alpha_{0} \sigma_{5} \alpha_{0}^{-1} \sigma_{4}^{-1}  \sigma_{5}
\sigma_{3}^{-1} \sigma_{4} \sigma_{2}^{-1} \alpha_{0}= \sigma_{2} \alpha_{0} \sigma_{5}^{-1} \alpha_{0} \sigma_{5}  \sigma_{3}^{-1}  \sigma_{4} \sigma_{2}^{-1} \sigma_{3}
\sigma_{1}^{-1}\\
&= \sigma_{2} \sigma_{1} \sigma_{2} \sigma_{3} \sigma_{4} \sigma_{1} \sigma_{2} \sigma_{3} \sigma_{4}
\sigma_{5}^2 \sigma_{4} \sigma_{3} \sigma_{3}^{-1} \sigma_{4}^{-1} \sigma_{3}^{-1}  \sigma_{4} \sigma_{2}^{-1} \sigma_{3} \sigma_{1}^{-1}\\
&= \sigma_{1} \sigma_{2} \sigma_{1} \sigma_{3} \sigma_{4} \sigma_{1}^{-1} \sigma_{2}^{-1} \sigma_{4}^{-1} \sigma_{3}^{-1}  \sigma_{2}^{-1} \sigma_{3}
\sigma_{1}^{-1}=1.
\end{align*}
This proves the claim, so $\ang{\gamma,\delta}\cong T_{1}$.
\end{rem}

\subsection{The geometric realisation of the finite subgroups of $B_{n}(\St)$}

The geometric realisation of the finite subgroups may be obtained by letting the corresponding subgroup of $\mathcal{M}_{0,n}$ act 
on the sphere with the $n$ strings attached in an appropriate manner. For the
subgroups $\dic{4n}$, $\Z_{2(n-1)}$ and $\dic{4(n-2)}$, we attach strings to $n$ symmetrically-distributed points (resp.\
$n-1$, $n-2$ points) on the equator, and $0$ (resp.\ $1$, $2$) points at the poles. For $T_1$, $O_1$ and $I$, the $n$
strings are attached symmetrically with respect to the associated regular polyhedron (for the values of $n$ given by
\reth{finitebn}) in the following manner. 

\item Let $H=\an[4]$ be the group of orientation-preserving symmetries of the tetrahedron. Then $n=6k+4$, $k\geq 0$, and
we take $k$ equally-spaced points in the interior of each edge, plus one point at each vertex (or face).

\item Let $H=\sn[4]$ be the group of orientation-preserving symmetries of the cube (or octahedron).
\begin{enumerate}[(i)]
\item $n=12k$, $k\in \N$: take $k$ equally-spaced points in the interior of each edge.
\item $n=12k+2$, $k\in \N$: take $k-1$ equally-spaced points in the interior of each edge, plus one point at each vertex and on each face.
\item $n=12k+6$, $k\geq 0$: take $k$ equally-spaced points in the interior of each edge, plus one point on each face.
\item $n=12k+8$, $k\geq 0$: take $k$ equally-spaced points in the interior of each edge, plus one point at each vertex.
\end{enumerate}
\item Let $H=\an[5]$ be the group of orientation-preserving symmetries of the icosahedron (or dodecahedron), which has
$12$ faces, $30$ edges and $20$ vertices.
\begin{enumerate}[(i)]
\item $n=30k$, $k\in \N$: take $k$ equally-spaced points in the interior of each edge.
\item $n=30k+2$, $k\in \N$: take $k-1$ equally-spaced points in the interior of each edge, plus one point at each vertex and on each face.
\item $n=30k+12$, $k\geq 0$: take $k$ equally-spaced points in the interior of each edge, plus one point on each face.
\item $n=30k+20$, $k\geq 0$: take $k$ equally-spaced points in the interior of each edge, plus one point at each vertex.
\end{enumerate}
\end{enumerate}

In each case, the action of the given group $H$ of symmetries yields the
corresponding maximal finite subgroup of $B_n(\St)$. This follows essentially
from the definition of the boundary operator
$\map{\partial}{\pi_1(D_n(\St))}[\pi_0\left(\diff{\St,X} \right)]$ in the long
exact sequence~\reqref{longexact} which we now describe in detail in our
setting. As in \resec{intro}, let $X$ be the basepoint in $D_n(\St)$, and let
$\map{\psi}{\diff{\St}}[D_n(\St)]$ denote evaluation on $X$. So if $g\in
\diff{\St}$ then $\psi(g)=g(X)$. Let $\id_{\St}$ be the basepoint in
$\diff{\St}$, so that $\psi(\id_{\St})=X$. Let $\beta \in B_n(\St)$ be a braid,
and let $\map{f}{[0,1]}[D_n(\St)]$ be a geometric braid which represents
$\beta$. So $f(0)=f(1)=X$, and the loop class $\ang{f}$ in $B_n(\St)$ is equal
to $\beta$. Then $f$ lifts to $\map{\widetilde{f}}{[0,1]}[\diff{\St}]$ which
satisfies $\widetilde{f}(0)=\id_{\St}$ and $\psi\circ \widetilde{f}=f$. Hence
$\psi\circ \widetilde{f}(1)=f(1)=X$, and thus $\widetilde{f}(1)$ belongs to the
fibre $\diff{\St,X}$. Geometrically, $\widetilde{f}$ is an isotopy of $\St$
which realises $\beta$ on the points of $X$. Neither $\widetilde{f}$ nor the
corresponding endpoint $\widetilde{f}(1)$ are unique, however all of the
possible $\widetilde{f}(1)$ belong to the same connected component of
$\diff{\St,X}$, and so determine a unique element, denoted $[\widetilde{f}(1)]$,
of $\pi_0\left(\diff{\St,X} \right)$, which is the image under $\partial$ of
$\beta$. Thus if $\widetilde{f}$ is an isotopy of $\St$ which realises $\beta$,
$\partial(\beta)$ is the mapping class of the homeomorphism $\widetilde{f}(1)$,
and corresponds geometrically to just remembering the final homeomorphism (in
particular, one forgets the strings of $\beta$).

Conversely, if $g\in \diff{\St}$ satisfies $g(X)=X$, let
$\map{h}{[0,1]}[\diff{\St}]$ be an isotopy from $h(0)=\id_{\St}$ to $h(1)=g$.
Then $\psi\circ h$ is a loop in $D_n(\St)$ based at $X$, so describes a
geometric braid obtained by attaching strings at the points of $X$ and following
the isotopy $h$. In $\St\times [0,1]$, the strings are given by $\brak{(\psi
\circ h(t),t)}_{t\in [0,1]}= \brak{(h(t)(X),t)}_{t\in [0,1]}$. Thus
$\ang{\psi\circ h}\in B_n(\St)$ is a braid, and by the above construction,
$\partial(\ang{\psi\circ h})=[h(1)]=[g]$. In other words, a choice of isotopy
$h$ between the identity and $g\in \diff{\St,X}$ allows us to lift the mapping
class $[g]$ to a preimage $\beta=\ang{\psi\circ h}$ under $\partial$ which is
obtained geometrically by attaching strings to $X$ during the isotopy $h$. 

Let $\map{r}{[0,1]}[\diff{\St}]$ denote rigid rotation through an angle $2\pi$. So $r(0)=r(1)=\id_{\St}$, the loop class
$\ang{r}$ generates $\pi_1\left(\diff{\St}\right) \cong \Z_2$, and thus $\ang{\psi\circ r}=\psi_{\ast}(\ang{r})= \ft$
since $\map{\psi_{\ast}}{\pi_1\left(\diff{\St}\right)}[B_n(\St)]$ is injective. The second preimage of $[g]$ under
$\partial$ is obtained by considering the isotopy $\map{h'}{[0,1]}[\diff{\St}]$ which is the isotopy $h$ followed by
$r$. The braids $\ang{\psi\circ h}$ and $\ang{\psi\circ h'}$ differ by $\ang{\psi\circ r}= \ft$, and thus define the two
preimages of $[g]$ under $\partial$.

Finally, each finite subgroup $H$ of $\mathcal{M}_{0,n}$ is realised by a finite subgroup of isometries of $\St$ (which
are the finite subgroups of $\operatorname{SO}(3)$)~\cite{K}. Each element of $H$ admits two preimages in $B_n(\St)$
which differ by $\ft$. These preimages thus make up the finite subgroup $\partial^{-1}(H)$ of $B_n(\St)$ whose order is
twice that of $H$.

\section{Position of the finite subgroups of $B_{n}(\St)$ relative to Murasugi's classification}\label{sec:muragi}

Let $n\geq 4$ be even. For $i=0,1,2$, let $G_i$ be the set of torsion
elements of $B_n(\St)$ whose order divides $2(n-i)$. Equivalently, by
\reth{murasugi}, $G_i$ is the set of conjugates of powers of
$\alpha_i$. Notice that $G_i$ is invariant under conjugation, $G_i\cap
G_j=\ang{\ft}$ for all $0\leq i<j\leq 2$, and $G_0\cup G_1 \cup G_2$
is the set of torsion elements of $B_n(\St)$. For many purposes, it is often useful to know
where a finite subgroup $H$ of $B_{n}(\St)$ lies relative to the $G_{i}$. In this section, we carry out
this calculation for all such subgroups.

\begin{prop}\label{prop:classgi}
Let $H$ be a finite subgroup of $B_{n}(\St)$ of order greater than or equal to $3$.
\begin{enumerate}[(I)] 
\item Suppose that $H$ is cyclic.
\begin{enumerate}[(a)]
\item if $\lvert H \rvert =4$ and $n$ is even then there exists a subgroup $H'$ of $B_{n}(\St)$ isomorphic to $\Z_{4}$
non-conjugate to $H$. One of $H,H'$ lies in $G_{0}$, while the other lies in $G_{2}$.
\item if either $\lvert H \rvert =4$ and $n$ is odd, or if $\lvert H \rvert \neq 4$ then $H
\subset G_{i}$, where $\lvert H \rvert \divides 2(n-i)$, and $i\in \brak{0,1,2}$.
\end{enumerate}
\item Suppose that $H$ is a subgroup of a maximal non-cyclic subgroup of $B_{n}(\St)$.
\begin{enumerate}[(a)]
\item If $H$ is a non-cyclic subgroup contained in $\dic{4n}$ or $\dic{4(n-2)}$ then it is itself dicyclic, of the form
$\dic{4k}$, where $k>1$ divides $n$ or $n-2$ respectively. Further:
\begin{enumerate}[(i)]
\item if $n$ is odd then $H\subset G_{i}\cup G_{1}$, where $i\in \brak{0,2}$ and $\lvert H\rvert
\divides 4(n-i)$.
\item Suppose that $n$ is even.
\begin{enumerate}[(1)]
\item if $k \divides n$ (resp.\ $k \divides n-2$) but $k\ndivides \frac{n}{2}$ (resp.\ $k\ndivides
\frac{n-2}{2}$) then $H$ lies in $G_{0}\cup G_{2}$ and meets both $G_{0}$ and $G_{2}$.
\item if $k\divides \frac{n}{2}$ (resp.\ $k\divides \frac{n-2}{2}$) then there exists another
subgroup $H'$ of $B_{n}(\St)$ isomorphic to $\dic{4k}$ but non conjugate to $H$. In this case, one
of $H,H'$ is contained wholly within $G_{0}$ (resp.\ $G_{2}$), and the other lies in $G_{0}\cup
G_{2}$ and meets both $G_{0}$ and $G_{2}$.
\end{enumerate}
\end{enumerate}
\item Suppose that $H$ is a subgroup of a copy of $T_{1}$ in the case that $T_{1}$ is maximal. 
\begin{enumerate}[(i)]
\item If $H\cong T_{1}$ then $H$ lies in $G_{0}\cup G_{1}$ (resp.\ $G_{2}\cup G_{1}$) if $n\equiv
4\bmod 12$ (resp.\ $n\equiv 10 \bmod 12$), and meets both $G_{0}$ (resp.\ $G_{2}$) and $G_{1}$.
\item If $H$ is isomorphic to $\Z_{3}$ or $\Z_{6}$ then it is contained in $G_{1}$.
\item If $H$ is isomorphic to $\Z_{4}$ or $\quat$ then it is contained in $G_{0}$ if $n\equiv 4\bmod
12$, and in $G_{2}$ if $n\equiv 10\bmod 12$.
\end{enumerate}
\item Suppose that $H$ is a subgroup of a copy of $I$ in the case that $I$ is maximal. 
\begin{enumerate}[(i)]
\item If $H$ is isomorphic to $I$ then $H$ is contained in $G_{0}$ (resp.\ $G_{2}$) if
$n\equiv 0 \bmod 60$ (resp.\ $n\equiv 2 \bmod 60$), and lies in $G_{0}\cup G_{2}$ and meets both
$G_{0}$ and $G_{2}$ if $n\equiv 12,20, 30,32,42,50 \bmod 60$.
\item If $H$ is isomorphic to $\Z_{3}$ or $\Z_{6}$ then it is contained in $G_{0}$ if $n\equiv
0, 12 \bmod 30$, and in $G_{2}$ if $n\equiv 2,20 \bmod 30$.
\item If $H$ is isomorphic to $\Z_{5}$ or $\Z_{10}$ then it is contained in
$G_{0}$ if $n\equiv 0, 20 \bmod 30$, and in $G_{2}$ if $n\equiv 2,12 \bmod 30$.
\item If $H$ is isomorphic to $\Z_{4}$ or $\quat$ then it is contained in $G_{0}$ if $n\equiv
0,12,20,32 \bmod 60$, and in $G_{2}$ if $n\equiv
2,30,42,50 \bmod 60$. 
\item If $H$ is isomorphic to $T_{1}$ or to $\dic{12}$ then it lies in $G_{0}$ if $n\equiv 0,12 \bmod
60$, in $G_{2}$ if $n\equiv 2,50 \bmod 60$, and lies in $G_{0}\cup G_{2}$ and meets both $G_{0}$ and
$G_{2}$ if $n\equiv 20, 30,32,42 \bmod 60$. 
\item If $H$ is isomorphic to $\dic{20}$ then it lies in $G_{0}$ if $n\equiv 0,20
\bmod 60$, in $G_{2}$ if $n\equiv 2,42 \bmod 60$, and lies in $G_{0}\cup G_{2}$ and meets both
$G_{0}$ and $G_{2}$ if $n\equiv 12,30,32,50 \bmod 60$.
\end{enumerate}
\item Suppose that $H$ is a subgroup of a copy of $O_{1}$ in the case that $O_{1}$ is maximal. 
\begin{enumerate}[(i)]
\item If $H$ is isomorphic to $O_{1}$ then it lies in $G_{0}$ if $n\equiv 0 \bmod 24$, in $G_{2}$ if
$n\equiv 2 \bmod 24$, and lies in $G_{0}\cup G_{2}$ and meets both $G_{0}$ and $G_{2}$ if $n\equiv
6,8,12,14, 18,20 \bmod 24$.
\item If $H$ is isomorphic to $T_{1}$ then it lies in $G_{0}$ if $n\equiv 0 \bmod 12$, in $G_{2}$
if $n\equiv 2 \bmod 12$, and lies in $G_{0}\cup G_{2}$ and meets both $G_{0}$ and $G_{2}$ if $n\equiv
6,8 \bmod 12$.
\item If $H$ is isomorphic to $\quat[16]$ then it lies in $G_{0}$ if $n\equiv 0,8 \bmod 24$, in
$G_{2}$ if $n\equiv 2,18 \bmod 24$, and lies in $G_{0}\cup G_{2}$ and meets both $G_{0}$ and $G_{2}$
if $n\equiv 6,12,14,20 \bmod 24$.
\item If $H$ is isomorphic to $\dic{12}$ then it lies in $G_{0}$ if $n\equiv 0,6 \bmod 24$, in
$G_{2}$ if $n\equiv 2,20 \bmod 24$, and lies in $G_{0}\cup G_{2}$ and meets both $G_{0}$ and $G_{2}$
if $n\equiv 8,12,14, 18 \bmod 24$.
\item If $H$ is isomorphic to $\Z_{8}$ then it lies in $G_{0}$ if $n\equiv 0,8 \bmod 12$, and in
$G_{2}$ if $n\equiv 2,6 \bmod 12$.
\item If $H$ is isomorphic to $\Z_{4}$ then there exists another non-conjugate subgroup $H'$ of
$B_{n}(\St)$ isomorphic to $\Z_{4}$. One of $H,H'$ is contained in $G_{0}$ if $n\equiv 0,8
\bmod{12}$, and in $G_{2}$ if $n\equiv 2,6 \bmod{12}$, while the other is contained in $G_{0}$ if
$n\equiv 0,6,8,14 \bmod 24$, and to $G_{2}$ if $n\equiv 2,12,18,20 \bmod 24$.
\item If $H$ is isomorphic to $\quat$ then there exists another non-conjugate subgroup $H'$ of
$B_{n}(\St)$ isomorphic to $\quat$. One of $H,H'$ is contained in $G_{0}$ if $n\equiv 0,8
\bmod{12}$, and to $G_{2}$ if $n\equiv 2,6 \bmod{12}$, while the other lies in $G_{0}$
if $n\equiv 0,8 \bmod 24$, in $G_{2}$ if $n\equiv 2,18 \bmod 24$, and lies in $G_{0}\cup G_{2}$ and
meets both $G_{0}$ and $G_{2}$ if $n\equiv 6,12,14,20 \bmod 24$.
\item If $H$ is isomorphic to $\Z_{3}$ or $\Z_{6}$ then it lies in $G_{0}$ if $n\equiv 0 \bmod 6$ and
in $G_{2}$ if $n\equiv 2 \bmod 6$.
\end{enumerate}
\end{enumerate}
\end{enumerate}
\end{prop}

\begin{proof}
Let $H$ be a finite subgroup of $B_{n}(\St)$ of order at least three. 
\begin{enumerate}[(I)]
\item Suppose first that $H$ is cyclic. 
Since $G_i\cap G_j=\ang{\ft}$ and $\lvert \ang{\alpha_{i}} \rvert =2(n-i)$, the order of $H$
is sufficient to decide where $H$ lies, unless $n$ is even and $H$ is of order $4$, in which case
there is another non-conjugate subgroup $H'$ isomorphic to $\Z_{4}$. One of 
$H,H'$ is conjugate to $\ang{\alpha_{0}^{n/2}}$ which is contained in $G_{0}$, while the other is
conjugate to  $\ang{\alpha_{2}^{(n-2)/2}}$ which lies in $G_{2}$. These two cases may be
distinguished easily by checking the permutation of a generator of $H,H'$.

\item Now suppose that $H$ is a subgroup of a maximal non-cyclic subgroup of $B_{n}(\St)$. We
consider the possible cases in turn.
\begin{enumerate}[(a)]
\item Firstly, let $H$ be a subgroup of the dicyclic group $\dic{4n}$, which up to conjugation may be
assumed to be $\ang{\alpha_{0}, \garside}= \ang{\alpha_{0}} \coprod \garside \ang{\alpha_{0}}$. We
first suppose that $n$ is odd. Then $\ang{\alpha_{0}} \subset G_{0}$, and the coset $\garside
\ang{\alpha_{0}}$ consists of the elements of $\dic{4n}$ of order $4$, so lies in $G_{1}$. The group
$\dic{4n}$ fits into a short exact sequence:
\begin{equation*}
1 \to \Z_{n} \to \dic{4n} \stackrel{g}{\to} \Z_{4} \to 1.
\end{equation*}
If $g(H)=\brak{\overline{0}}$, then $H< \Z_{n}$, and $H$ is cyclic, of order dividing $n$, so lies
in $G_{0}$. If $g(H)=\brak{\overline{0}, \overline{2}}$, then $H< \Z_{2n}$, and again $H$ is cyclic,
of order dividing $2n$, so lies in $G_{0}$. Finally, if $g(H)=\Z_{4}$ then we have
\begin{equation*}
1 \to H\cap \Z_{n} \to H \stackrel{g}{\to} \Z_{4} \to 1,
\end{equation*}
and $H\cong \Z_{k} \rtimes \Z_{4}$, where $k$ divides $n$. If $k=1$ then $H\cong \Z_{4}$. Since $n$
is odd, $H$ must then lie in $G_{1}$. So suppose that $k>1$. Then $H=\ang{\alpha_{0}^{n/k},
\garside}$ is dicyclic, and so lies in $G_{0}\cup G_{1}$.

Now suppose that $n$ is even. Then $\dic{4n}$ fits into the following short exact sequence:
\begin{equation*}
1 \to \Z_{2n} \to \dic{4n} \stackrel{f}{\to} \Z_{2} \to 1.
\end{equation*}
If $f(H)=\brak{\overline{0}}$ then
$H\subset \Z_{2n}$ and so lies in $G_{0}$. If $f(H)=\Z_{2}$ and $H\cap \Z_{2n}$ were of odd order, then $H$ would be
both dicyclic and of order twice an odd number, which cannot occur. So suppose that $f(H)=\Z_{2}$ and $H\cap \Z_{2n}$
is of even order, $2k$, say, where $k\divides n$. If $k=1$ then $H\cong \Z_{4}$, and $H$ may lie in $G_{0}$ or $G_{2}$
depending on the permutation of its generators. So suppose that $k\geq 2$. Then $H$ is dicyclic of order $4k$. Now 
\begin{equation*}
\dic{4n} = \underbrace{\ang{\alpha_{0}}}_{\subset G_{0}} \coprod
\underbrace{\garside \ang{\alpha_{0}^2}}_{\subset G_{0}} \coprod
\underbrace{\garside \alpha_{0}\ang{\alpha_{0}^2}}_{\subset G_{2}}.
\end{equation*}
The inclusions follow from the fact that the elements of $\garside \ang{\alpha_{0}^2}$ (resp.\
$\garside \alpha_{0}\ang{\alpha_{0}^2}$) are conjugate (in $\dic{4n}$), $\garside \in G_{0}$, and 
\begin{align*}
\pi(\garside \alpha_{0})&= (1 , n) (2 , n-1) \cdots \left(\frac{n}{2} , \frac{n}{2}+1 \right) (1
, n , \ldots ,2)\\
&= (n)\left(\frac{n}{2} \right)(1,n-1) (2 , n-2) (3 , n-3) \cdots \left(\frac{n}{2}-1 , \frac{n}{2}+1
\right),
\end{align*}
where $\map{\pi}{B_{n}(\St)}[\sn]$ denotes the homomorphism defined  on the generators by $\pi(\sigma_{i})= (i, i+1)$.
Thus $\garside \alpha_{0} \in G_{2}$.

If $k \ndivides \frac{n}{2}$ then by \repr{isoconj}, there is just one conjugacy class of $\dic{4k}$ of the form
$\ang{\alpha_{0}^{n/k}, \garside}$, and since $n/k$ is odd, we have
\begin{equation*}
\dic{4k} = \underbrace{\ang{\alpha_{0}^{n/k}}}_{\subset G_{0}} \coprod
\underbrace{\garside \ang{\alpha_{0}^{n/k}}}_{\subset G_{2}}.
\end{equation*}
In particular, all of the elements of $\dic{4k}$ of order $4$ belong to $G_{2}$. Thus $\dic{4k} \cap
(G_{0}\setminus G_{2}) \neq \vide$ and $\dic{4k} \cap (G_{2}\setminus G_{0}) \neq \vide$.

If $k \divides \frac{n}{2}$ then by \repr{isoconj}, there are two non-conjugate copies of $\dic{4k}$ given
by
\begin{equation*}
\ang{\alpha_{0}^{n/k}, \garside}= \underbrace{\ang{\alpha_{0}^{n/k}}}_{\subset G_{0}} \coprod
\underbrace{\garside \ang{\alpha_{0}^{n/k}}}_{\subset G_{0}}, \;\text{and} 
\end{equation*}
and 
\begin{equation*}
\ang{\alpha_{0}^{n/k}, \garside\alpha_{0}}=\underbrace{\ang{\alpha_{0}^{n/k}}}_{\subset G_{0}}
\coprod \underbrace{\garside \alpha_{0}\ang{\alpha_{0}^{n/k}}}_{\subset G_{2}}.
\end{equation*}
The first copy lies entirely within $G_{0}$, while the second lies in $G_{0}\cup G_{2}$ and meets both $G_{0}\setminus
G_{2}$ and $G_{2} \setminus G_{0}$. 

A similar result holds for $\dic{4(n-2)}$: its subgroups are either subgroups of $\Z_{2(n-2)}$,
so lie in $G_{2}$, or else are dicyclic, of the form $\dic{4k}$, where $k \divides n-2$. If $k=1$
then the subgroup in question is $\ang{\garside}$ which lies in $G_{0}$. If $k>1$ then as above,
we distinguish two cases. If $k \ndivides \frac{n-2}{2}$ then there is just one copy of $\dic{4k}$
which lies in $G_{0}\cup G_{2}$ and meets both $G_{0}\setminus G_{2}$ and $G_{2} \setminus G_{0}$. If $k \divides
\frac{n-2}{2}$, then setting $\alpha_{2}' =\alpha_{0} \alpha_{2} \alpha_{0}^{-1}$,
there are two copies of $\dic{4k}$, $\ang{\alpha_{2}'^{n/k}, \garside}$, which lies in $G_{0}\cup G_{2}$ and meets both
$G_{0}\setminus G_{2}$ and $G_{2} \setminus G_{0}$, and $\ang{\alpha_{2}'^{n/k}, \alpha_{2}'\garside}$, which is
contained in $G_{2}$. 

\item Suppose that $H$ is a subgroup of a copy of
$T_{1}$ when $T_{1}$ is maximal, so $n \equiv 4 \bmod 6$. Assume first that $H\cong T_{1}$. Since
$H\cong \quat \rtimes \Z_{3}$, all of its order $4$ elements are conjugate, and so all elements of
$\quat$ must lie in the same $G_{i}$. Now $\quat=\dic{8}$, so from above, we must be in one of the
cases $2 \divides \frac{n}{2}$ or $2 \divides \frac{n-2}{2}$. Indeed if $n\equiv 4 \bmod 12$ then
$n=4+12l=4(1+3l)$, $l\in \N$, and so $\quat$ is contained in $G_{0}$, while if $n\equiv 10 \bmod 12$
then $n=10+12l=2(5+6l)$, $l\in \N$, and so $\quat$ is contained in $G_{2}$. The remaining elements of
$H$ are of order $3$ or $6$, and since $n\equiv 4 \bmod 6$, lie in $G_{1}$. So if $n\equiv 4\bmod
12$ (resp.\ $n\equiv 10 \bmod 12$) then $H$ lies in $G_{0}\cup G_{1}$ (resp.\ $G_{2}\cup G_{1}$)
and meets both $G_{0}$ (resp.\ $G_{2}$) and $G_{1}$.

From this, we deduce immediately the following: if $H$ is isomorphic to $\Z_{3}$ or $\Z_{6}$ then
it is contained in $G_{1}$, and if it is isomorphic to $\Z_{4}$ or $\quat$ then it is contained in
$G_{0}$ if $n\equiv 4\bmod 12$, and in $G_{2}$ if $n\equiv 10\bmod 12$.

\item Suppose that $H$ is a subgroup of a copy of
$I$ when $I$ is maximal, so $n\equiv 0,2,12,20 \bmod
30$. Assume first that $H\cong I$. So $I$ has a subgroup isomorphic to $T_{1}$, whose copy of
$\quat$ lies entirely in $G_{0}$ or $G_{2}$. The subgroups of order $8$ of $H$ are its Sylow
$2$-subgroups, so are conjugate, and thus all lie either in $G_{0}$ or in $G_{2}$. Hence from the
analysis of the dicyclic case, $2$ divides $\frac{n}{2}$ or $\frac{n-2}{2}$. Further, all
elements of $H$ of order $4$ are contained in one of its subgroups isomorphic to $\quat$ (because the order $2$ elements
of $\an[5]$ are the product of two transpositions, and are contained in a subgroup isomorphic to $\Z_{2} \oplus \Z_{2}$, which lifts to $\quat$ in $I$).
Hence all order $4$ elements of $H$ lie either in $G_{0}$ if $4 \divides n$, or in $G_{2}$ if $4
\divides n-2$. The remaining elements of $H$ are of order $3,6,5$ and $10$, and lie in either $G_{0}$
or $G_{2}$ depending on the value of $n$ modulo the order. Thus $H$ lies
entirely in $G_{0}$ (resp.\ $G_{2}$) if $n\equiv 0 \bmod 60$ (resp.\ $n\equiv 2 \bmod 60$), and lies
in $G_{0}\cup G_{2}$ and meets both $G_{0}$ and $G_{2}$ if $n\equiv 12,20, 30,32,42,50 \bmod 60$.
 
We now consider the other possibilities for subgroups of $I$: if $H$ is isomorphic to either $\Z_{3}$
or $\Z_{6}$, it is contained in $G_{0}$ if $n\equiv 0, 12 \bmod 30$, and in $G_{2}$ if $n\equiv 2,20
\bmod 30$;  if $H$ is isomorphic to either $\Z_{5}$ or $\Z_{10}$, it is contained in
$G_{0}$ if $n\equiv 0, 20 \bmod 30$, and in $G_{2}$ if $n\equiv 2,12 \bmod 30$; and if $H$ is
isomorphic to either $\Z_{4}$ or $\quat$, it is contained in $G_{0}$ if $n\equiv 0,12,20,32 \bmod 60$,
and in $G_{2}$ if $n\equiv 2,30,42,50 \bmod 60$. Next, if $H$ is isomorphic to $T_{1}$, it
consists of a copy of $\quat$ and elements of order $3$ and $6$, so lies in $G_{0}$ if $n\equiv 0,12
\bmod 60$, in $G_{2}$ if $n\equiv 2,50 \bmod 60$, and lies in $G_{0}\cup G_{2}$ and meets both
$G_{0}$ and $G_{2}$ if $n\equiv 20, 30,32,42 \bmod 60$. Now suppose that $H$ is isomorphic to
$\dic{12}\cong \Z_{3} \rtimes \Z_{4}= \Z_{6} \coprod \garside \Z_{6}$. Since the elements
of $\garside \Z_{6}$ are of order $4$, it follows from the analysis of the cyclic subgroups that
$H$ satisfies the same conditions as in the case of $T_{1}$. Finally, if $H$ is isomorphic to
$\dic{20}\cong \Z_{5} \rtimes \Z_{4}= \Z_{10} \coprod \garside \Z_{10}$, since the elements of
$\garside \Z_{10}$ are of order $4$, it follows from the analysis of the cyclic subgroups that
$H$ lies in $G_{0}$ if $n\equiv 0,20 \bmod 60$, in $G_{2}$ if $n\equiv 2,42 \bmod 60$, and lies in
$G_{0}\cup G_{2}$ and meets both $G_{0}$ and $G_{2}$ if $n\equiv 12,30,32,50 \bmod 60$.

\item Suppose that $H$ is a subgroup of a copy of $O_{1}$ when $O_{1}$ is maximal, so
$n\equiv 0,2 \bmod 6$. Assume first that $H\cong O_{1}$. Then it has a subgroup isomorphic to $T_{1}$
(which is unique since $\sn[4]$ has a unique subgroup abstractly isomorphic to $\an[4]$), and the
copy of $\quat$ lying in $T_{1}$ lies entirely in $G_{0}$ if $n\equiv 0,8 \bmod 12$, and in $G_{2}$
if $n\equiv 2,6 \bmod 12$. The complement of this copy of $\quat$ in $T_{1}$ consists of elements of
order $3$ and $6$, and so lie in $G_{0}$ if $n\equiv 0 \bmod 6$ and in $G_{2}$ if $n\equiv 2 \bmod 6$
(thus the subgroups of $O_{1}$ isomorphic to $\Z_{3}$ and $\Z_{6}$ lie in $G_{0}$ if $n\equiv 0 \bmod
6$ and in $G_{2}$ if $n\equiv 2 \bmod 6$). Thus $T_{1}$ lies in $G_{0}$ if $n\equiv 0 \bmod 12$, in
$G_{2}$ if $n\equiv 2 \bmod 12$, and lies in $G_{0}\cup G_{2}$ and meets both $G_{0}$ and $G_{2}$ if
$n\equiv 6,8 \bmod 12$.

In order to analyse the remaining possible subgroups $\quat[16]$, $\dic{12}$, $\dic{20}$ of $O_{1}$,
as well as the other copy of $\quat$ lying in $\quat[16]$, we must study the elements of $H \setminus
T_{1}$. They project to elements of $\sn[4] \setminus \an[4]$, which are either $4$-cycles,
or transpositions. We analyse the geometric formulation of $O_{1}$ described in \resec{geomreal} as being
obtained from the action of $\sn[4]$ on a cube, with the $n$ strings attached appropriately. The
$4$-cycles are realised by rotations by $\pi/2$ about an axis which passes through the centres of
two opposite faces. This gives rise to an element of $G_{0}$ if the $n$ marked points are not
these central points (i.e.\ if $n\equiv 0,8 \bmod 12$), and to elements of $G_{2}$ if some 
of the $n$ marked points are central points of the faces (i.e.\ if $n\equiv 2,6 \bmod 12$). The
transpositions are realised by rotations by $\pi$ about an axis which passes through the centres of
two diagonally-opposite edges. This gives rise to an element of $G_{0}$ if there are an even number of marked
points on each edge (i.e.\ if $n\equiv 0,6,8,14 \bmod 24$), and to elements of $G_{2}$ if there
are an odd number of marked points on each edge (i.e.\ if $n\equiv 2,12,18,20 \bmod 24$). Putting
together these results with those for $T_{1}$, if $H\cong O_{1}$, we conclude that it lies in $G_{0}$ if $n\equiv 0
\bmod 24$, in $G_{2}$ if $n\equiv 2 \bmod 24$, and lies in $G_{0}\cup G_{2}$ and meets both $G_{0}$
and $G_{2}$ if $n\equiv 6,8,12,14, 18,20 \bmod 24$.

Now suppose that $H$ is a subgroup of a copy of $O_{1}$ isomorphic to $\quat[16]$. Such subgroups are
the Sylow $2$-subgroups of $O_{1}$, so are conjugate. If $n\equiv 0 \bmod 24$  (resp.\ $n\equiv 2
\bmod 24$) then $O_{1}$ lies in $G_{0}$ (resp.\ $G_{2}$), and hence so does $\quat[16]$. So suppose
that $n\nequiv 0,2 \bmod 24$. Any subgroup of $O_{1}$ isomorphic to $\quat[16]$ contains elements of
order $8$ which lie in $O_{1}\setminus T_{1}$, and so are associated with the above $4$-cycles.
Further, $H$ projects to a subgroup of $\sn[4]$ isomorphic to $\dih{8}$ which is
generated by a $4$-cycle and a transposition. Studying the associated rotations as above, if one has
fixed points and the other not then automatically $H$ lies in $G_{0}\cup G_{2}$ and meets
both $G_{0}$ and $G_{2}$. This occurs when $n\equiv 6,12,14,20 \bmod 24$. So suppose that $n\equiv
8,18 \bmod 24$.

If $n\equiv 8 \bmod 24$ (resp.\ $n\equiv 18 \bmod 24$) then the elements of $H$ corresponding
to the $4$-cycles and the transpositions of $\dih{8}$ belong to $G_{0}$ (resp.\ $G_{2}$). Further,
the remaining elements of $\dih{8}$ are products of such elements, and so the corresponding elements
in $H$ are also elements of $T_{1}\cong \quat \rtimes \Z_{3}$ of order $4$. But such elements
lie in the $\quat$-factor. Since $n\equiv 8 \bmod 12$ (resp.\ $n \equiv 6 \bmod 12$), this copy of
$\quat$ lies in $G_{0}$ (resp.\ $G_{2}$), and hence so does the given subgroup $\quat[16]$.
Summing up, $H$ lies in $G_{0}$ if $n\equiv 0,8 \bmod 24$, in $G_{2}$ if $n\equiv 2,18 \bmod
24$, and lies in $G_{0}\cup G_{2}$ and meets both $G_{0}$ and $G_{2}$ if $n\equiv 6,12,14,20
\bmod 24$.

Now suppose that $H$ is a subgroup of a copy of $O_{1}$ isomorphic to $\dic{12}$. If $n\equiv 0
\bmod 24$  (resp.\ $n\equiv 2 \bmod 24$) then $O_{1}$ lies in $G_{0}$ (resp.\ $G_{2}$), and hence so
does $H$. So suppose that $n\nequiv 0,2 \bmod 24$. Any subgroup of $O_{1}$ isomorphic to
$H$ projects onto a subgroup of $\sn[4]$ isomorphic to $\sn[3]$ which consists of $3$-cycles and
transpositions. Hence $H$ is generated by an element of order $4$ lying in $O_{1} \setminus
T_{1}$, and an element of order $6$, which lies in $T_{1}$. The first element belongs to $G_{0}$ if
$n\equiv 6,8,14 \bmod 24$ and to $G_{2}$ if $n\equiv 12,18,20 \bmod 24$, while the second element
belongs to $G_{0}$ if $n\equiv 6,12,18 \bmod 24$ and to $G_{2}$ if $n\equiv 8,14,20 \bmod 24$.
Hence if $n\equiv 8,12,14,18 \equiv 24$ then $H$ lies in $G_{0}\cup G_{2}$ and meets both
$G_{0}$ and $G_{2}$. The product of the two given generators is also of order $4$ and so lies in
$G_{0}$ if $n\equiv 6 \bmod 24$, and in $G_{2}$ if $n\equiv 20 \bmod 24$. Thus $H$ lies in
$G_{0}$ if $n\equiv 0,6 \bmod 24$, in $G_{2}$ if $n\equiv 2,20 \bmod 24$, and lies in $G_{0}\cup
G_{2}$ and meets both $G_{0}$ and $G_{2}$ if $n\equiv 8,12,14, 18 \bmod 24$.

Now suppose that $H$ is a subgroup of a copy of $O_{1}$ isomorphic to $\Z_{4}$. There are two
possibilities. If it is contained in the copy of $\quat$ lying in the subgroup $T_{1}$, from the
results for $\quat$, we see that $H$ lies in $G_{0}$ if $n\equiv 0,8 \bmod{12}$, and in
$G_{2}$ if $n\equiv 2,6 \bmod{12}$. The second possibility is that $H$ possesses elements in
$O_{1}\setminus T_{1}$, and emanates from the rotation of order $2$ whose permutation is a
transposition. Thus it is contained in $G_{0}$ if $n\equiv 0,6,8,14 \bmod 24$, and to $G_{2}$ if $n\equiv
2,12,18,20 \bmod 24$. 

Finally, suppose that $H$ is a subgroup of a copy of $O_{1}$ isomorphic to $\quat$. Again there 
are two possibilities. If $H$ lies in the subgroup $T_{1}$, it is contained in $G_{0}$ if $n\equiv
0,8 \bmod{12}$, and to $G_{2}$ if $n\equiv 2,6 \bmod{12}$. The second possibility is that
it projects to a subgroup of $\sn[4]$ generated by two transpositions having disjoint support.
Such a subgroup thus has four elements of order $4$ in $O_{1}\setminus T_{1}$ and two in $T_{1}$.
From the results obtained in the case of $\Z_{4}$, we see that $H$ lies in $G_{0}$
if $n\equiv 0,8 \bmod 24$, in $G_{2}$ if $n\equiv 2,18 \bmod 24$, and lies in $G_{0}\cup G_{2}$ and
meets both $G_{0}$ and $G_{2}$ if $n\equiv 6,12,14,20 \bmod 24$.
\end{enumerate}
\end{enumerate}
\end{proof}

\section{Realisation of finite groups as subgroups of the lower central and derived series of $B_{n}(\St)$}\label{sec:gamma2}

In this section, we consider the realisation of the finite subgroups of \reth{finitebn} as subgroups of the lower
central $\Gamma_i(B_n(\St))$ and derived series $(B_n(\St)^{(i)})$ of $B_{n}(\St)$. By~\cite{GG6}, we already know that
if $4\mid n$ then $\Gamma_2(B_n(\St))$ has a subgroup isomorphic to $\quat$. If $n\geq 4$ is even but not divisible by
$4$, we may ask if the same result is true if $4 \nmid n$. We start by proving \reth{gammaconj}, which is the case of the
dicyclic groups. We then then complete the analysis of the other finite subgroups in \repr{gamma2}. 

\begin{proof}[Proof of \reth{gammaconj}] 
Suppose that $n$ is even. Let $N\in \brak{n-2,n}$, set $N=2^lk$ where $l\in \N$ and $k$ is odd, and let $x=\alpha_0$
(resp.\ $x=\alpha_0 \alpha_2 \alpha_0^{-1}$) if $N=n$ (resp.\ $N=n-2$).

\begin{enumerate}[(a)]
\item Since $B_n(\St)$ has a subgroup $\ang{x,\garside}$ isomorphic to
$\dic{4N}=\dic{2^{l+2}k}$, the statement is true for $j=0$.  So suppose the
result holds for some $j\in \brak{0,1,\ldots, l-1}$. Then $B_n(\St)$ contains
$2^j$ copies of $\dic{2^{l+2-j}k}$ of the form $\ang{x^{2^j}, x^i \garside}$,
for $i=0,1, \ldots, 2^j-1$. Hence $\ang{x^{2^{j+1}}, x^i\garside}$ is a subgroup of $\ang{x^{2^j}, x^i \garside}$ isomorphic to $\dic{2^{l+1-j}k}$. But since
\begin{gather*}
\left( x^{(2^j+i)} \garside \right)^2 = x^{(2^j+i)} \garside x^{(2^j+i)} \garside^{-1} \garside^2=\ft, \quad\text{and}\\
x^{(2^j+i)} \garside \cdot x^{2^{j+1}} \left( x^{(2^j+i)} \garside \right)^{-1} = x^{-2^{j+1}},
\end{gather*}
it follows that $\ang{x^{2^{j+1}}, x^{(2^j+i)} \garside}$ is also a subgroup of $\ang{x^{2^j}, x^i \garside}$ isomorphic to $\dic{2^{l+1-j}k}$.

If $q$ is any divisor of $k$, then replacing $x$ by $x^q$ yields also $2^j$ copies $\ang{x^{2^jq}, x^{iq}\garside}$,
$i=0,1,\ldots, 2^j-1$, of $\dic{2^{l+2-j}k/q}$ for $j\in \brak{0,1,\ldots, l}$.

\item If $j=0$, the statement holds trivially. So suppose that $j\geq 1$. From part~(\ref{it:2jcopies}), $\ang{x^{2^jq},
x^{iq}\garside}$ and $\ang{x^{2^jq}, x^{i'q}\garside}$ are subgroups of $B_n(\St)$ isomorphic to $\dic{2^{l+2-j}k/q}$.
Under the Abelianisation homomorphism $\map{\xi}{B_n(\St)}[\Z_{2(n-1)}]$, $\xi(x)=  \overline{n-1}$, and
\begin{equation*}
\xi(\garside)=\xi\bigl((\sigma_1\cdots \sigma_{n-1})\cdots (\sigma_1 \sigma_2)\sigma_1\bigr)= \overline{\frac{1}{2} n(n-1)}= 
\begin{cases}
\overline{0} & \text{if $\frac{n}{2}$ is even}\\
\overline{n-1} & \text{if $\frac{n}{2}$ is odd.}
\end{cases}
\end{equation*}
Since $j\geq 1$, $\xi(x^{2^jq})= \overline{0}$. Furthermore,
 \begin{equation*}
\xi \left( x^{iq}\garside \right)=
\begin{cases}
\overline{0} & \text{if $\frac{n}{2}+i$ is even}\\
\overline{n-1} & \text{if $\frac{n}{2}+i$ is odd.}
\end{cases}
\end{equation*}
So $\ang{x^{2^jq}, x^{iq}\garside}\subset \Gamma_2(B_n(\St))$ if and only if $\frac{n}{2}+i$ is even. Thus if $i-i'$ is
odd, the subgroups $\ang{x^{2^jq}, x^{iq}\garside}$ and $\ang{x^{2^jq}, x^{i'q}\garside}$ cannot be conjugate. But by
\repr{isoconj}(\ref{it:conj2}), these are precisely the conjugacy classes of subgroups isomorphic to
$\dic{2^{l+2-j}k/q}$. The result follows.
\end{enumerate}
\end{proof}

From this, we may deduce \repr{quatgamma2}.

\begin{proof}[Proof of \repr{quatgamma2}]
We use the notation of the proof of \reth{gammaconj}. If $j=0$ and $q$ is an odd divisor of $n$ then there is just one
conjugacy class of the abstract group $\dic{4n/q}$, which is realised as $\ang{x^q,\garside}$. Now $x^q\notin
\Gamma_2(B_n(\St))$, so $\dic{4n/q} \nsubset \Gamma_2(B_n(\St))$.

If $j\geq 1$ then as we saw in the proof of \reth{gammaconj}, $\ang{x^{2^jq}, x^{iq}\garside}\subset \Gamma_2(B_n(\St))$
if and only if $\frac{n}{2}+i$ is even. So with $i=0,1$, one of $\ang{x^{2^jq}, \garside}$
and $\ang{x^{2^jq}, x^q\garside}$ is contained in $\Gamma_2(B_n(\St))$, while the other is not. 

Finally, let $N$ be the element of $\brak{n,n-2}$ divisible by $4$. Then $l\geq 2$, and taking $q=k$ and $j=l-1$, from
the previous paragraph, one of  $\ang{x^{N/2}, \garside}$ and $\ang{x^{N/2}, x^k\garside}$ (the two non-conjugate copies of $\quat$)
belongs to $\Gamma_2(B_n(\St))$, the other not.
\end{proof}

We now give the analogous result for the cyclic and binary polyhedral subgroups of $B_{n}(\St)$.
\begin{prop}\label{prop:gamma2}
Let $G$ be a finite subgroup of $B_{n}(\St)$. 
\begin{enumerate}[(a)]
\item Suppose that $G$ is cyclic. 
\begin{enumerate}[(i)]
\item If $G$ is of order~$2$, then $G\subset \Gamma_{2}(B_{n}(\St))$ if and only if $n$ is even.
\item Suppose that $G$ is of order greater than or equal to $3$. Then either:
\begin{enumerate}[--]
\item $\lvert G \rvert$ divides $2(n-1)$ in which case $G\nsubset \Gamma_{2}(B_{n}(\St))$, or
\item $\lvert G \rvert$ divides $2(n-i)$, where $i\in\brak{0,2}$. In this case, $G \subset \Gamma_{2}(B_{n}(\St))$ if and only if $\lvert G
\rvert$ divides $n-i$.
\end{enumerate}
\end{enumerate}
\item Suppose that $G$ is a subgroup of order at least $3$ of some binary polyhedral subgroup  $H$ of $B_{n}(\St)$.
\begin{enumerate}[(i)]
\item Suppose that $H \cong T_{1}$ in the case that $T_{1}$ is maximal. Then $G\subset
\Gamma_{2}(B_{n}(\St))$ if $G\cong \Z_{4},\quat$, and $G\nsubset \Gamma_{2}(B_{n}(\St))$ if $G\cong \Z_{3}, \Z_{6},
T_{1}$.
\item Suppose that $H \cong I$ in the case that $I$ is maximal. Then $G\subset
\Gamma_{2}(B_{n}(\St))$.
\item Suppose that $H \cong O_{1}$ in the case that $O_{1}$ is maximal. If $G$ is
contained in the subgroup $K$ of $H$ isomorphic to $T_{1}$ then $G\subset \Gamma_{2}(B_{n}(\St))$. If $G\nsubset K$
then $G\subset \Gamma_{2}(B_{n}(\St))$ if $n\equiv 0,2,8,18 \bmod 24$, and $G\nsubset \Gamma_{2}(B_{n}(\St))$ if $n\equiv 6,12,14,20 \bmod 24$.
\end{enumerate}
\end{enumerate}
\end{prop}

\begin{proof}
We set $\Gamma_{2}=\Gamma_{2}(B_{n}(\St))$. If $G$ is of order~$2$, then $G=\ang{\ft}$ and as $\xi(\ft)=\overline{n(n-1)}$, it follows easily that $G\subset
\Gamma_{2}$ if and only if $n$ is even. We assume from now on that $\lvert G \rvert \geq 3$. Since $\Gamma_{2}$ is
normal in $B_{n}(\St)$, we may work up to conjugation.

First suppose that $G$ is cyclic. Then by \reth{murasugi}, it is conjugate to a subgroup of $\ang{\alpha_{i}}$ for some
$i \in \brak{0,1,2}$. If $i=1$ then $\xi(\alpha_{1}^{j})=\overline{jn}$ for all $j\in \Z$. If $\alpha_{1}^{j}\in
\Gamma_{2}$ then there exists $k\in \Z$ such that $jn=2k(n-1)$, thus $n-1 \divides j$, and so $j=l(n-1)$ for some
$l\in\Z$. But then $\alpha_{1}^j=\alpha_{1}^{l(n-1)} \in \ang{\ft}$. We conclude that $\ang{\alpha_{1}}\cap \Gamma_{2}
\subset \ang{\ft}$. Hence $G \nsubset \Gamma_{2}$.

Suppose then that $G$ is conjugate to a subgroup of $\ang{\alpha_{i}}$, where $i=0,2$. Set $k=\lvert G \rvert$. Then $\xi(\alpha_{i})=\overline{n-1}$,
$k\divides 2(n-i)$, and up to conjugacy, $G=\ang{\alpha_{i}^{2(n-i)/k}}$. So $G \subset \Gamma_{2}$
if and only if $2(n-i)/k$ is even, which is equivalent to $k \divides n-i$. Thus if $G$ is conjugate to a subgroup of $\ang{\alpha_{i}}$, where
$i=0,2$, we have:
\begin{equation}\label{eq:gamma2}
G \subset \Gamma_{2} \Longleftrightarrow \lvert G \rvert  \divides n-i.
\end{equation}

Now suppose that $H$ is isomorphic to $T_{1}$ in the case that $T_{1}$ is maximal, so that $n\equiv 4
\bmod 6$. If $G$ is isomorphic to $T_{1}, \Z_{6}$ or $\Z_{3}$ then the order $3$ elements lie in $G_{1}\setminus
\ang{\ft}$, and from the cyclic case, it follows that $G \nsubset \Gamma_{2}$. So assume that $G$ is isomorphic to either $\Z_{4}$
or $\quat$. Since $\quat$ is generated by elements of order $4$, it suffices to analyse the case $\Z_{4}$. By
\repr{classgi}, $G$ lies in $G_{0}$ if $n \equiv 4 \bmod 12$, and in $G_{2}$ if $n \equiv 10 \bmod 12$. In both cases,
$G\subset \Gamma_{2}$ by \req{gamma2}.

Now suppose that $H$ is isomorphic to $I$ in the case that $I$ is maximal, so that $n\equiv 0,2,12,20
\bmod 30$. We claim that $G\subset \Gamma_{2}$ whatever the value of $n$. To see this, it suffices to check that all
of the maximal cyclic subgroups $\Z_{4}$, $\Z_{6}$, $\Z_{10}$ of $I$ are contained in $\Gamma_{2}$. This follows
easily from \repr{classgi} and \req{gamma2}. 
 
Now suppose that $H$ is isomorphic to $O_{1}$ in the case that $O_{1}$ is maximal, so that $n\equiv
0,2 \bmod 6$. Again it suffices to consider the maximal cyclic subgroups $\Z_{4}$, $\Z_{6}$ and $\Z_{8}$ of $O_{1}$.
Applying \repr{classgi} and \req{gamma2}, we obtain the following results:
\begin{enumerate}[--]
\item if $G$ is isomorphic to $\Z_{8}$, it projects to a subgroup of $\sn[4]$ generated by a $4$-cycle. Then $G\subset
G_{0}$ if $n\equiv 0,8 \bmod 12$, and $G\subset G_{2}$ if $n\equiv 2,6 \bmod 12$, and so $G\subset \Gamma_{2}$ if $n\equiv
0,2,8,18 \bmod 24$, and $G\nsubset \Gamma_{2}$ if $n\equiv 6,12,14,20 \bmod 24$.
\item if $G$ is isomorphic to $\Z_{6}$ then $G\subset \Gamma_{2}$.
\item if $G$ is isomorphic to $\Z_{4}$, there are two possibilities. If $G$ lies in the subgroup $K$ of $O_{1}$
isomorphic to $T_{1}$ then $G\subset \Gamma_{2}$. Otherwise $G$ is generated by an element of order $4$ not belonging to
$K$, in which case we obtain the same answer as for $\Z_{8}$.
\end{enumerate}
Since every cyclic subgroup of order $3$ of $O_{1}$ is contained in one of order $6$, this gives the results if $G$ is
cyclic. Suppose now that $G=K$. Then $G$ is generated by the elements of order $6$ and the elements of order $4$
belonging to $K$, so $G\subset \Gamma_{2}$. 

If $G$ is abstractly isomorphic to $\quat[16]$ then it is generated by elements of order $8$, elements of order $4$
lying in $K$, and elements of order $4$ not lying in $K$. From above, we have that $G\subset \Gamma_{2}$ if $n\equiv
0,2,8,18 \bmod 24$, and $G\nsubset \Gamma_{2}$ if $n\equiv 6,12,14,20 \bmod 24$.

If $G$ is abstractly isomorphic to $\quat$ then there are two possibilities: either $G$ lies in $K$, so is contained
in $\Gamma_{2}$, or else it is generated by elements of order $4$ not belonging to $K$. In this case, from above,
$G\subset \Gamma_{2}$ if $n\equiv 0,2,8,18 \bmod 24$, and $G\nsubset \Gamma_{2}$ if $n\equiv 6,12,14,20 \bmod 24$.

Finally, suppose that $G$ is abstractly isomorphic to $\dic{12}$. Then it projects to a copy of $\sn[3]$ in $\sn[4]$.
From above, it follows that $G\subset \Gamma_{2}$ if $n\equiv 0,2,8,18 \bmod 24$, and $G\nsubset \Gamma_{2}$ if $n\equiv 6,12,14,20 \bmod 24$.
 \end{proof}
 
\begin{rem}
Having dealt with the behaviour of the finite subgroups relative to the commutator subgroup of $B_{n}(\St)$, one might
ask what happens for the higher elements of the lower central series $\brak{\Gamma_{i}(B_{n}(\St))}_{{n\in\N}}$ and of the
derived series $\brak{\left( B_{n}(\St) \right)^{(i)}}_{i\geq 0}$ of $B_{n}(\St)$. But if $n\neq 2$ (resp.\ $n\geq
5$), the lower central series (resp.\ derived series) of $B_{n}(\St)$ is stationary from the commutator subgroup
onwards~\cite{GG5}. It just remains to look at the derived series of $B_{4}(\St)$. Recall from that paper that
$(B_{4}(\St))^{(1)}$ is a semi-direct product of $\quat$ by a free group $\F[2]$ of rank two, that $(B_{4}(\St))^{(2)}$
is a semi-direct product of $\quat$ by the derived subgroup $(\F[2])^{(1)}$ of $\F[2]$, that $(B_{4}(\St))^{(3)}$
is the direct product of $\ang{\ft[4]}$ by $(\F[2])^{(2)}$, and that $(B_{4}(\St))^{(i+1)} \cong (\F[2])^{(i)}$ for
all $i\geq 3$. Thus there is a copy of $\quat$ which lies in $(B_{4}(\St))^{(2)}$ but not in $(B_{4}(\St))^{(3)}$. The
full twist remains until $(B_{4}(\St))^{(3)}$, and then $(B_{4}(\St))^{(4)}$ is torsion free.
\end{rem}

\section{A proof of Murasugi's theorem}\label{sec:murasugi}

Let $H_1,H_2$ be isomorphic finite cyclic subgroups of
$\mathcal{M}_{0,n}$. From \reth{stukow}, if $n$ is odd, or if $n$ is
even and $\lvert H_1 \rvert= \lvert H_2 \rvert\neq 2$ then $H_1$ and
$H_2$ are conjugate. If $n$ is even, there are exactly two conjugacy
classes of subgroups of $\mathcal{M}_{0,n}$ of order $2$, and thus
there are exactly two conjugacy classes of subgroups of $B_n(\St)$ of
order $4$.

It follows from \resec{class} that:

\begin{prop}\label{prop:isoconjbn}
Let $G_1,G_2$ be isomorphic finite cyclic subgroups of order $m$ of $B_n(\St)$. If $n$ is odd, or if $n$ is even and
$m\neq 4$ then $G_1$ and $G_2$ are conjugate. If $n$ is even, there are exactly two conjugacy classes of subgroups of
$B_n(\St)$ of order $4$.\qed
\end{prop}

If $n$ is even then $\alpha_0^{n/2}$ and $\alpha_2^{(n-2)/2}$ are of order $4$, and they generate non-conjugate
subgroups since their images in $\sn$ are not conjugate, which yields the two conjugacy classes of $\Z_4$ of
\repr{isoconjbn}. From this, we may deduce \reth{murasugi}.

\begin{proof}[Proof of \reth{murasugi}]
Let $x\in B_n(\St)$ be a torsion element. Then $\ang{x}$ is contained in a maximal cyclic subgroup $C$ of one of the
maximal finite subgroups $G$ of $B_n(\St)$ given by \reth{finitebn}.

First suppose that $n$ is odd. Then $G$ is one of $\Z_{2(n-1)}$, $\dic{4n}$ and $\dic{4(n-2)}$, and so $C$ must be one
of $\Z_{2(n-1)}$, $\Z_{2n}$, $\Z_{2(n-2)}$ and $\Z_4$. Hence $C$ is isomorphic to $\ang{\alpha_1}$, $\ang{\alpha_0}$,
$\ang{\alpha_2}$ and $\ang{\alpha_1^{(n-1)/2}}$ respectively. So by \repr{isoconjbn}, $x$ is conjugate to a power of one
of $\alpha_0$, $\alpha_1$ and $\alpha_2$ which proves the theorem in this case.

Now suppose that $n$ is even. If $C \cong \Z_4$ then $C$ is conjugate to one of $\ang{\alpha_0^{n/2}}$ or
$\ang{\alpha_2^{(n-2)/2}}$, and the result holds. So suppose that $C \ncong \Z_4$. If $G$ is one of $\Z_{2(n-1)}$,
$\dic{4n}$ and $\dic{4(n-2)}$, then $C$ is one of $\Z_{2(n-1)}$, $\Z_{2n}$, $\Z_{2(n-2)}$, and so is isomorphic to
$\ang{\alpha_1}$, $\ang{\alpha_0}$, and $\ang{\alpha_2}$ respectively. If $G=T_1$ (so $n\equiv 4 \bmod 6$) then $C \cong
\Z_6$, and so is conjugate to $\ang{\alpha_1^{(n-1)/3}}$. If $G=O_1$ (so $n\equiv 0,2 \bmod 6$) then $C \cong \Z_6$ or
$C \cong \Z_8$, and so is conjugate to $\ang{\alpha_0^{n/3}}$ or $\ang{\alpha_2^{(n-2)/3}}$. Finally, if $G=I$ (so
$n\equiv 0,2,12,20 \bmod 30$) then $C$ is isomorphic to one of $\Z_6$ or $\Z_{10}$. If $C\cong \Z_6$ then $C$ is
conjugate to $\ang{\alpha_0^{n/3}}$ if $n\equiv 0,12 \bmod 30$ or to $\ang{\alpha_2^{(n-2)/3}}$ if $n\equiv 2,20 \bmod
30$. If $C\cong \Z_{10}$ then $C$ is conjugate to $\ang{\alpha_0^{n/5}}$ if $n\equiv 0,20 \bmod 30$ or to
$\ang{\alpha_2^{(n-2)/5}}$ if $n\equiv 2,12 \bmod 30$. In all cases, $x$ is conjugate to a power of one of $\alpha_0$,
$\alpha_1$ and $\alpha_2$, which completes the proof of the theorem.
\end{proof}

\end{document}